# RIGID ANALYTIC FLATIFICATORS

Hans Schoutens

ABSTRACT. Let $K$ be an algebraically closed field endowed with a complete non-archimedean norm. Let $f\colon Y \to X$ be a map of $K$-affinoid varieties. We prove that for each point $x \in X$, either $f$ is flat at $x$, or there exists, at least locally around $x$, a maximal locally closed analytic subvariety $Z \subset X$ containing $x$, such that the base change $f^{-1}(Z) \to Z$ is flat at $x$, and, moreover, $g^{-1}(Z)$ has again this property in any point of the fibre of $x$ after base change over an arbitrary map $g\colon X' \to X$ of affinoid varieties. If we take the local blowing up $\pi\colon \tilde{X} \to X$ with this centre $Z$, then the fibre with respect to the strict transform $\tilde{f}$ of $f$ under $\pi$, of any point of $\tilde{X}$ lying above $x$, has grown strictly smaller. Among the corollaries to these results we quote, that flatness in rigid analytic geometry is local in the source; that flatness over a reduced quasi-compact rigid analytic variety can be tested by surjective families; that an inclusion of affinoid domains is flat in a point, if it is unramified in that point.

## 0. INTRODUCTION

The goal of this article is to prove the existence of universal flatificators in the setting of rigid analytic geometry and to investigate the behaviour of fibres when we blow up with centre such a universal flatificator. Let us be more precise.

We fix once and for all an algebraically closed field $K$ endowed with a complete non-archimedean norm. The analytic-geometric objects we will be interested in are the *rigid analytic varieties*. These are locally ringed spaces which look locally like an *affinoid variety* $\mathrm{Sp}\, A$, the latter being the collection of all maximal ideals of an affinoid algebra $A$ on which $A$ acts as a function algebra. Recall that an *affinoid algebra* is a quotient of some free Tate algebra $K\langle X \rangle$, where the latter is defined as the collection of all formal power series over $K$ in the variables $X$ for which the general coefficient tends to zero. For more details, we refer to **[3]**.

We define a *rigid analytic flatificator* of a map $f\colon Y \to X$ of rigid analytic varieties at a point $x \in X$ for which the fibre $f^{-1}(x) \neq \emptyset$, as a locally closed immersion $i\colon Z \hookrightarrow X$ containing $x$ for which the base-change $Y \times_X Z \to Z$ becomes flat and such that if $j\colon V \hookrightarrow X$ is any locally closed immersion containing $x$ for which $Y \times_X V \to V$ is flat at $x$, then $j$ factors through $i$ at $x$, i.e., $V \hookrightarrow Z$ at $x$. Here we will say that a map is flat at a point of the target space, if it is flat in each point of the fibre over this point.

We call such a flatificator *universal* at $x$, if, for any map $X' \to X$ of rigid analytic varieties, the locally closed immersion $X' \times_X Z \hookrightarrow X'$ obtained from $i$ by base change, is a rigid analytic flatificator of the base change $X' \times_X Y \to X'$ of $f$ at any point of $X'$ lying above $x$. Our first main theorem guarantees the existence of such a universal flatificator at each point, for an arbitrary map of quasi-compact rigid analytic varieties.

**Theorem.** *Let $f\colon Y \to X$ be a map of quasi-compact rigid analytic varieties and let $x \in X$ be an arbitrary point of $X$ with non-empty fibre. Then there exists an admissible affinoid $U$ of $X$ containing $x$, such that the restriction $f^{-1}(U) \to U$ has a universal rigid analytic flatificator at $x$.*

*Proof.* See **(3.7)**.

As a corollary to this theorem, we show that flatness is an open condition in the source: given a map $Y \to X$ of rigid analytic varieties which is flat in a point $y$ of $Y$, then there exists an admissible affinoid $U$ of $Y$ containing $y$, such that the restriction $U \to X$ is flat. See **(3.8)**. A second corollary asserts that one can use surjective families to test flatness over reduced quasi-compact rigid analytic varieties: let $\pi_i\colon X_i \to X$ be maps of rigid analytic spaces, with $X$ reduced and quasi-compact, such that the images $\{\mathrm{Im}(\pi_i)\}$ of the $\pi_i$ cover $X$ (point set wise), then any map $f\colon Y \to X$ of quasi-compact rigid analytic varieties is flat, if and only if, all the base changes $Y \times_X X_i \to X_i$ are flat. See **(3.9)**.

The existence of flatificators in the complex analytic case had been proven already by Hironaka in [**9**], using Weierstrass Preparation. Our proof is based on his arguments, but uses a very general construction to prove the existence of a flatificator (reviewed in the first two sections), which has enabled us to prove similar results in an algebraic-geometric setting in [**17**]. In his paper, Hironaka further investigates the behaviour of a blowing up map with centre a universal flatificator. We obtain in the same fashion the analogous Fibre Lemma, which roughly says that if one blows up $X$ with centre a universal rigid analytic flatificator $Z$ at a point $x$ of a map $f\colon Y \to X$, then the fibre of $f$ over $x$ is strictly bigger than the corresponding fibre of the strict transform of $f$ over a point in the blowing up lying above $x$. Using that one can always modify the centre locally in such way that it becomes nowhere dense, without disturbing the fibres, we obtain the following theorem.

**Theorem.** *Let $f\colon Y \to X$ be a map of quasi-compact rigid analytic varieties and let $x \in X$ be an arbitrary point. If $X$ is reduced and $f$ is not flat at $x$, then there exists a local blowing up $\pi\colon \tilde{X} \to X$ with nowhere dense centre $Z$ containing $x$, such that the restriction $f^{-1}(Z) \to Z$ is flat and such that, if*

$$\begin{array}{ccc} \tilde{Y} & \xrightarrow{\theta} & Y \\ \tilde{f} \downarrow & & \downarrow f \\ \tilde{X} & \xrightarrow{\pi} & X \end{array}$$

*denotes the diagram of the strict transform of $f$ under $\pi$, then we have, for every $\tilde{x} \in \pi^{-1}(x)$, a non-trivial closed immersion of fibres*

$$\tilde{f}^{-1}(\tilde{x}) \subsetneq f^{-1}(x).$$

*Proof.* See **(5.3)**.

As a corollary we want to mention the following theorem (see **(5.5)**): given an inclusion $A \hookrightarrow B$ of affinoid domains and a maximal ideal $\mathfrak{m}$ of $A$ which extends to a maximal ideal $\mathfrak{m}B$ in $B$, then $A \to B$ is flat at $\mathfrak{m}$. In other words, in this situation, unramified is the same as etale.

However, the main application of this theorem will appear in a future paper [**7**] with Gardener on Quantifier Elimination for the rigid analytic language expanded with a truncated division function, where we will combine his work on the Voûte Etoilée in [**6**] and our present results on rigid analytic flatificators. Let me give a brief review of the contents of this paper in which we study the structure of rigid subanalytic sets. We reduce this to the study of images of affinoid maps. In the flat case, one knows, due to a result of Raynaud that they are semianalytic. Hence the goal is to *flatten* an arbitrary map of affinoid varieties $f\colon Y \to X$ by means of local blowing up maps. This is first done locally in some point $x \in X$ where $f$ is not flat: take its flatificator and use this as centre of a local blowing up. If the strict transform of $f$ is flat at any point above $x$, then we are done, otherwise continue. However, by our second main Theorem **(5.3)**, the fibre has grown strictly smaller and hence the tree of local blowing ups used to render $f$ flat has finite depth. Exploiting the compactness properties of the Voûte Etoilée one can then show that a finite tree already suffices and by another compactness argument the global case follows.

However, for all this to work we need to carry out the above flattening procedure in the Berkovich spaces associated to the affinoid varieties rather than in the affinoid varieties themselves. So we will need the Berkovich counterpart of the theory of rigid analytic flatificators. A brief outline of how various of our results translate to this setting is given in an appendix. Berkovich spaces were introduced by Berkovich in [**1**] and [**2**] to reestablish a genuine Hausdorff topology on a rigid analytic space. This characteristic was exploited by Gardener to define the aforementioned Voûte Etoilée. The appendix shows basically that all results carry over to this new setup modulo some minor modifications. In the proofs we give merely an outline of these modifications, more details can be found in Gardener's D.Phil. thesis [**5**].

*0.2. Remark.* A note on the base field $K$. Although we have taken $K$ to be algebraically closed, this is not necessary for any of the results of this paper. The reader can check that we have virtually nowhere used this assumption. We have used it in **(3.8)**, where we used that $K$ is equal to the residue field $A/\mathfrak{m}$, for $\mathfrak{m}$ a maximal ideal in an affinoid algebra, where in the general case this is a finite field extension. But since any extension of fields is faithfully flat, one can easily reduce to the algebraic closed case for this theorem. The only other difference upon dropping the algebraic closedness condition is a notational one. The affinoid space associated to a free Tate algebra $K\langle S\rangle$, with $S = (S_1,\ldots,S_N)$, can no longer be identified with the Cartesian product $R^N$, but instead one has to take the closed unit ball $\mathbb{B}^N = \mathbb{B}^N(\bar{K})$, modulo the action of the absolute Galois group $\operatorname{Gal}(\bar{K}/K)$, where $\bar{K}$ is the algebraic closure of $K$. This amounts in replacing $R^N$ or $R$ by $\mathbb{B}^N$ or $\mathbb{B}^1$, wherever the former occur, such as in **(3.4-6)**).

**0.3. Notation and Terminology.** All rigid analytic varieties will be over $K$. The valuation ring of $K$ will be denoted by $R$, which sometimes will be identified with the unit disk, i.e., with $\operatorname{Sp}(K\langle T\rangle)$, where $T$ is one variable. We call a separated rigid analytic variety admitting a finite admissible affinoid covering, *quasi-compact*.

We will freely interchange the terms *closed immersion* and *closed analytic subvariety* of a rigid analytic variety $X$. Any such $Z$ is uniquely defined by and uniquely defines a coherent $\mathcal{O}_X$-ideal. If we are only interested in the underlying point set of $Z$, then we will write $|Z|$ and refer to the latter as a *closed analytic subset* of $X$. When $Z$ and $Z'$ are closed analytic subvarieties of $X$, defined respectively by $\mathcal{I}$ and $\mathcal{I}'$, then $Z \cap Z'$ stands for the closed analytic subvariety $Z \times_X Z'$, i.e., defined by the sum $\mathcal{I} + \mathcal{I}'$. Likewise, if $f\colon Y \to X$ is a map of rigid analytic varieties, then $f^{-1}(Z)$ stands for the closed analytic variety $Y \times_X Z$, i.e., defined by $\mathcal{I}\mathcal{O}_Y$. Whenever, however, we take unions, direct images or complements (except for the notion of an *analytic complement*, introduced in Section 4), we will disregard any analytic structure and consider only the underlying point sets. A locally closed analytic subvariety of a rigid analytic variety $X$ is a closed analytic subvariety of an admissible affinoid open of $X$. Similarly defined are locally closed analytic subsets.

Blowing up maps are defined in full generality in rigid analytic geometry; for a reference, see our paper [**18**] . By a *local blowing up map* we mean the composition of a blowing up map and an affinoid open immersion, i.e., if $U$ is an admissible affinoid of $X$ and $\theta\colon \tilde{X} \to U$ is a blowing up map, then the composition $\tilde{X} \to U \hookrightarrow X$ is a local blowing up map.

The topological notions of *closure*, *openness* and *nowhere density* are always taken with respect to the Zariski-topology.

**0.4. Acknowledgment.** We want to thank Tim Gardener for the useful and stimulating conversations on the subject.

1. GENERAL FLATIFICATORS

**1.1. Definition.** Let us briefly recall the definition of the flatificator of a module over a local ring, as studied in [**17**] . In what follows, let $(S,\mathfrak{p})$ be a Noetherian local ring and let $M$ be an $S$-module. The *flatness filter* $\mathbb{F}\mathrm{l}_S(M)$ of $M$ over $S$ is defined as the collection of all ideals $I$ of $S$, such that $M/IM$ becomes a flat $S/I$-module. When this flatness filter happens to be principal, i.e., has a least element, say $\mathfrak{f}$, then we call $\mathfrak{f}$ the *flatificator* of $M$ over $S$.

From [**17**, Proposition 2.3] , it follows that, if $M$ has a founded and flat $\mathfrak{p}$-representation, then it has a flatificator. Here, for an arbitrary (not necessarily local) ring $A$, an $A$-module $F$ is called *founded*, if for any submodule $H$ of $F$, we have that $H \subset \mathfrak{f}_H F$, where $\mathfrak{f}_H$ is the intersection of all ideals $I$ of $A$, for which $H \subset IF$. An $A$-module $M$ is said to have a *blue $\mathfrak{p}$-representation*, where blue stands for any property of modules and where $\mathfrak{p}$ is a prime ideal of $A$, if there exists an exact sequence of $A$-modules

$$(1) \qquad 0 \to H \longrightarrow F \longrightarrow M \to 0,$$

such that $F$ is blue and $H \subset \mathfrak{p}F$. In , it has been proven that if (1) is such a founded and flat $\mathfrak{p}$-representation for $S$ and $M$ as above, then $\mathfrak{f}_H$ is the flatificator of $M$.

We will be mainly concerned in this paper with modules over an affinoid algebra. Let us therefore introduce the following definition. Let $X = \operatorname{Sp} A$ be an affinoid variety and let $\mathcal{F}$ be an $\mathcal{O}_X$-module. We will say that $\mathcal{F}$ is *strongly founded*, if $\mathcal{F}(X)$ is a founded $A$-module and, moreover, for all $x \in X$, the module $\mathcal{F}_x$ is founded as an $\mathcal{O}_{X,x}$-module. It is not too hard to prove that a finite direct sum of (strongly) founded modules is again such a module.

**1.2. Proposition.** *Let $A$ be an affinoid algebra and let $T = (T_1, \ldots, T_n)$ be a finite number of variables. Then $A\langle T \rangle$ is a strongly founded $A$-module.*

*Proof.* Let us first show that $A\langle T \rangle$ is founded. Let $H$ be an $A$-submodule of $A\langle T \rangle$ and define $\mathfrak{f}_H$ as above. Let $f = \sum_i a_i T^i \in H$. Then for any ideal $I$ of $A$ such that $H \subset IA\langle T \rangle$, we have that $f \in IA\langle T \rangle$ and hence each $a_i \in I$. Therefore each $a_i \in \mathfrak{f}_H$ and hence $f \in \mathfrak{f}_H A\langle T \rangle$.

We now show that the same holds true in the stalk over a point $x$ of $X = \operatorname{Sp} A$. Let $\mathcal{F}$ denote the $\mathcal{O}_X$-module associated to $A\langle T \rangle$. One easily verifies that $\mathcal{F}_x = A\langle T \rangle \otimes_A \mathcal{O}_{X,x}$. Let $H \subset \mathcal{F}_x$ be an $\mathcal{O}_{X,x}$-submodule and, as always, let $\mathfrak{f} = \mathfrak{f}_H$ be the intersection of all $\mathcal{O}_{X,x}$-ideals $I$ such that $H \subset I\mathcal{F}_x$. Let $I$ be such an ideal and let $h \in H$. After passing to a small enough affinoid subdomain $U = \operatorname{Sp} C \subset X$ containing $x$ over which $I$ and $h$ are still defined, we can write

$$(1) \qquad h = \sum_{i<d} \sum_{\alpha \in \mathbb{N}^n} h_{i\alpha} T^\alpha \otimes s_i,$$

in $A\langle T \rangle \otimes_A C$, with $h_{i\alpha} \in A$ and $s_i \in C$, for $\alpha \in \mathbb{N}^n$ and $i < d$. Via the canonical injection $\mathcal{F}(U) = A\langle T \rangle \otimes_A C \hookrightarrow C\langle T \rangle$, we obtain that

$$\sum_{i,\alpha} h_{i\alpha} s_i T^\alpha \in IC\langle T \rangle,$$

and whence that $\sum_i h_{i\alpha} s_i \in I$, for all $\alpha$. Since this holds for all ideals $I$ such that $H \subset I\mathcal{F}_x$, we must have that $\sum_i h_{i\alpha} s_i \in \mathfrak{f}$. After possibly further shrinking $U$, we may assume that also $\mathfrak{f}$ is defined over $U$, so that using (1) we see that $h = 0$ as an element of $(C/\mathfrak{f}C)\langle T \rangle$ whence also in $A\langle T \rangle \otimes_A C/\mathfrak{f}C$, since the latter ring is embedded in the former. In other words, that $h \in \mathfrak{f}\mathcal{F}_x$, as wanted. ∎

**1.3. Definition.** Let $S$ be a ring, $M$ an $S$-module and $\mathfrak{a}$ an ideal in $S$. We say that an $S$-submodule $N$ of $M$ is $\mathfrak{a}$-*closed* if

$$N \cap \mathfrak{a}M = \mathfrak{a}N.$$

This is equivalent with the statement that

$$0 \to N \longrightarrow M \longrightarrow C \to 0$$

remains exact after tensoring it up with $S/\mathfrak{a}$, that is,

$$0 \to N/\mathfrak{a}N \longrightarrow M/\mathfrak{a}M \longrightarrow C/\mathfrak{a}C \to 0$$

is exact.

We call a chain of $S$-submodules $M_i$ of $M$,

$$M_0 \subset M_1 \subset \cdots \subset M_s = M$$

an $\mathfrak{a}$-*series*, if each $M_{i-1}$ is $\mathfrak{a}$-closed in $M_i$, for $i = 1, \ldots, s$.

**1.4. Lemma.** *Let $S$ be a ring, $\mathfrak{a}$ an ideal of $S$ and $M$ an $S$-module. Let $N \subset W \subset M$ be $S$-submodules and suppose that*

  (1) *$N \subset M$ is $\mathfrak{a}$-closed,*
  (2) *$W/N \subset M/N$ is $\mathfrak{a}$-closed.*

*Then $W \subset M$ is $\mathfrak{a}$-closed.*

*Proof.* Let $q \in \mathfrak{a}M \cap W$. If $\bar{q}$ denotes the residue of $q$ in $M/N$, then we have that $\bar{q} \in \mathfrak{a}M/N \cap W/N$ and hence by (2) that $\bar{q} \in \mathfrak{a}W/N$. In other words, there exist $w \in \mathfrak{a}W$ and $n \in N$ with $q = w + n$. But then we have that $n \in \mathfrak{a}M \cap N$, so that by (1) $n \in \mathfrak{a}N$ and hence $q \in \mathfrak{a}W$. ∎

**1.5. Lemma.** *Let $S$ be a ring, $M$ a finitely generated $S$-module and let $\mathfrak{p}$ be a prime ideal of $S$ belonging to $\mathrm{Supp}(M)$. Then there exists a non-zero cyclic, $\mathfrak{p}$-closed submodule $N$ of $M$.*

*Proof.* Since $M_\mathfrak{p} \neq 0$, we get by Nakayama's Lemma that

$$\mathfrak{p} M_\mathfrak{p} \neq M_\mathfrak{p},$$

meaning that $\mathfrak{p} \in \mathrm{Supp}(M/\mathfrak{p}M)$. But then $\mathfrak{p}$ is a minimal prime in the support of $M/\mathfrak{p}M$ and hence an associated prime of $M/\mathfrak{p}M$. So there exists an $\overline{\omega} \in M/\mathfrak{p}M$ such that

$$\mathrm{Ann}_S(\overline{\omega}) = \mathfrak{p}.$$

Let $\omega \in M$ be a lifting of $\overline{\omega}$ and set $N = S\omega$. We claim that $N$ is $\mathfrak{p}$-closed. For let $\alpha \in N \cap \mathfrak{p}M$, then there exists an $r \in S$ such that $\alpha = r\omega$. From $r\omega \in \mathfrak{p}M$, we have that $r\overline{\omega} = 0$ and hence $r \in \mathrm{Ann}_S(\overline{\omega}) = \mathfrak{p}$, which implies that $\alpha \in \mathfrak{p}N$. ∎

**1.6. Proposition.** *Let $S$ be a Noetherian ring, $M$ a finitely generated $S$-module and $\mathfrak{p}$ a prime ideal of $S$. Then there exists a $\mathfrak{p}$-series of $M$,*

$$(1) \qquad 0 = N_0 \subset N_1 \subset \cdots \subset N_t \subset M,$$

*where for each $i = 1, \ldots, t$ we have that $N_i/N_{i-1}$ is a cyclic $S$-module and*

$$(M/N_t)_\mathfrak{p} = 0.$$

*Proof.* Set $N_0 = 0$. If $M_\mathfrak{p} = 0$, then $0 \subset M$ is such a series, hence we may assume that $M_\mathfrak{p} \neq 0$, that is, $\mathfrak{p} \in \mathrm{Supp}(M)$. By **(1.5)**, we can find a non-zero cyclic, $\mathfrak{p}$-closed submodule $N_1$ of $M$. Let $Q_1 = M/N_1$. Either $(Q_1)_\mathfrak{p} = 0$, and hence $0 \subset N_1 \subset M$ is a $\mathfrak{p}$-series as desired. In the other case, we have that $\mathfrak{p} \in \mathrm{Supp}(Q_1)$ and by the same reasoning as above, we can find a submodule $N_2$ of $M$, such that $N_1 \subsetneq N_2$ and $N_2/N_1$ is cyclic and $\mathfrak{p}$-closed in $Q_1$. From **(1.4)** we get then that $N_2$ is $\mathfrak{p}$-closed in $M$. Let $Q_2 = M/N_2$, then again, either $(Q_2)_\mathfrak{p} = 0$ and hence

$$0 \subsetneq N_1 \subsetneq N_2 \subset M,$$

is a $\mathfrak{p}$-series of the required type, or $\mathfrak{p} \in \mathrm{Supp}(Q_2)$ and then we continue as above. However, since $M$ is Noetherian, this process has to stop and we end up with a $\mathfrak{p}$-series as proclaimed. ∎

## 2. Analytically Projective Modules

**2.1. Definition.** Let $A$ be a normed ring and $M$ a normed $A$-module. If $M$ is a normed $A$-module which is complete, then we will call $M$ an *analytic $A$-module*.

If $M$ and $N$ are two analytic $A$-modules, then we will say that a morphism $\varphi\colon M \to N$ is a *morphism of analytic $A$-modules*, if $\varphi$ is a continuous and bounded $A$-module morphism. Here, a morphism $\varphi$ is called *bounded*, if there exists a real positive constant $C$, such that, for all $m \in M$, we have that $|\varphi(m)| \leq C \cdot |m|$. The morphism is called *strict*, if the quotient topology on $\varphi(M)$ coincides with the topology on $\varphi(M)$ inherited from the topology on $N$. Moreover, if $\varphi$ is a surjective strict morphism, then using [**3**, 1.1.9. Lemma 2], we can show that, for any bounded sequence $(\nu_k)_k$ in $N$, we can find a bounded sequence $(\mu_k)_k$ in $M$, such that $\varphi(\mu_k) = \nu_k$.

**2.2. Proposition.** *Let $A$ be a normed ring and $M$ and $N$ two finite analytic $A$-modules. If $M$ and $N$ are pseudo-Cartesian, then any continuous $A$-module morphism $\varphi\colon M \to N$ is bounded. Moreover, if $\varphi$ is surjective, then it is strict.*

*Remark.* Recall that a finite normed $A$-module $M$ is called *pseudo-Cartesian*, if there exists a finite system of generators $\{\mu_1, \ldots, \mu_s\}$ of $M$, such that any element $\mu \in M$, can be written as $\mu = a_1\mu_1 + \cdots + a_s\mu_s$, for some $a_i \in A$, such that

$$(\dagger) \qquad |\mu| = \max_{1 \leq i \leq s}\{|a_i| \cdot |\mu_i|\}.$$

Such a system of generators is then called a *pseudo-Cartesian* system of generators.

From [**3**, 5.2.7] it follows that any finitely generated module over an affinoid algebra $A$ is pseudo-Cartesian and analytic, where we consider $A$ as a normed ring by taking any residue norm on it. Note that any such norm defines the same topology on $A$, and hence, for our purposes, can be chosen arbitrarily. In particular, any morphism between finitely generated modules over an affinoid algebra is a morphism of analytic $A$-modules, which is strict whenever it is surjective.

*Proof.* Let $\{\mu_1, \ldots, \mu_s\}$ be a pseudo-Cartesian system of generators for $M$. Obviously, we can always throw out any $\mu_i = 0$, and hence we may assume that none of the $\mu_i$ are zero. Let $\mu \in M$ and let $a_i \in A$, be such that $\mu = a_1\mu_1 + \cdots + a_s\mu_s$ and such that (†) holds. Let $C$ be the maximum of all the $|\varphi(\mu_i)|/|\mu_i|$. Hence, we obtain that

$$|\varphi(\mu)| \leq \max_i |a_i| \cdot |\varphi(\mu_i)|$$
$$\leq C \cdot \max_i |a_i| \cdot |\mu_i| = C|\mu|,$$

where the last equality comes from (†). This proves the boundedness.

Assume now that $\varphi$ is moreover surjective. Let $\varepsilon > 0$. According to [**3**, 1.1.9. Lemma 2], we need to find a $\delta > 0$, such that we can find, for any $\mu \in M$ with $|\varphi(\mu)| \leq \delta$, an element $\kappa \in \ker\varphi$, such that $|\mu + \kappa| \leq \varepsilon$. Let $\{\nu_1, \ldots, \nu_t\}$ be a pseudo-Cartesian system of (non-zero) generators for $N$. Choose $\gamma_i \in M$ such that $\varphi(\gamma_i) = \nu_i$. Let $D$ be the maximum of all $|\gamma_i|/|\nu_i|$. We claim that $\delta = \varepsilon/D$ will do the job. Indeed, take any $\mu \in M$ with $|\varphi(\mu)| \leq \varepsilon/D$. Write $\varphi(\mu) = \sum_i b_i\nu_i$, so that $|\varphi(\mu)| = \max_i |b_i| \cdot |\nu_i|$. Let $\gamma = \sum_i b_i\gamma_i$. Then $\gamma - \mu \in \ker\varphi$. Hence we are done if we can show that $|\gamma| \leq \varepsilon$. But we have that

$$|\gamma| \leq \max_i |b_i| \cdot |\gamma_i|$$
$$\leq D \cdot \max_i |b_i| \cdot |\nu_i|$$
$$= D|\varphi(\mu)| \leq \varepsilon.$$

∎

**2.3. Definition.** Let $A$ be a normed ring and $P$ an analytic $A$-module. We will say that $P$ is *analytically projective*, if for any strict surjective morphism $\varepsilon\colon M \to N$ of analytic $A$-modules and for any morphism $\theta\colon P \to N$ of analytic $A$-modules, we can find a morphism $\sigma\colon P \to M$ of analytic $A$-modules, making the following diagram commutative

$$\begin{array}{ccccc} M & \xrightarrow{\varepsilon} & N & \longrightarrow & 0 \\ \sigma\uparrow & & \uparrow\theta & & \\ P & =\!\!=\!\!= & P. & & \end{array}$$

It is an easy exercise to prove that a finite direct sum is analytically projective, if and only if, each of its summands is such.

**2.4. Proposition.** *If $A$ is an affinoid algebra, then $A\langle T\rangle$ is an analytically projective $A$-module.*

*Proof.* Let $\varepsilon\colon M \twoheadrightarrow N$ be a strict surjective morphism of analytic $A$-modules and let $\theta\colon A\langle T\rangle \to N$ be an arbitrary morphism of analytic $A$-modules. Let $\nu_i = \theta(T^i)$. By boundedness and $|T^i| = 1$, we have that the $|\nu_i|$ are bounded. Since $\varepsilon$ is strict, we can find $\mu_i \in M$, such that $|\mu_i|$ are bounded and $\varphi(\mu_i) = \nu_i$ (see the final remark in **(2.1)**). Let us therefore define $\sigma\colon A\langle T\rangle \to M$ as follows. If $f = \sum_i a_i T^i \in A\langle T\rangle$, then we set

$$\sigma(f) \stackrel{\text{def}}{=} \sum_i a_i\mu_i,$$

where this sum does converge since $|a_i| \to 0$ and the $|\mu_i|$ are bounded. It is now a straightforward computation to verify the requirements of the statement for this morphism $\sigma$. ∎

**2.5. Proposition.** *Let $A$ be a normed ring and $\mathfrak{p}$ a prime ideal in $A$. Let*

$$0 \to M_1 \to M \xrightarrow{\theta} M_2 \to 0$$

*be a short exact sequence of analytic $A$-modules. If $M_1$ is $\mathfrak{p}$-closed in $M$, $\theta$ is strict and both $M_1$ and $M_2$ have an analytically projective $\mathfrak{p}$-representation, then also $M$ has an analytically projective $\mathfrak{p}$-representation.*

*Moreover, if the two former representations are blue, where blue is a property of $A$-modules which is closed under taking finite direct sums, then also the latter representation is blue.*

*Proof.* Let $\pi_i \colon F_i \twoheadrightarrow M_i$ be an analytically projective $\mathfrak{p}$-representation of $M_i$, for $i = 1, 2$. In other words, $F_i$ is an analytically projective $A$-module and $\ker(\pi_i) \subset \mathfrak{p} F_i$. Now, since $F_2$ is analytically projective and $\theta$ is strict, we can find a morphism $\psi_2$ making the following diagram commutative

$$\begin{array}{ccc} F_2 & \xrightarrow{\psi_2} & M \\ \| & & \downarrow \theta \\ F_2 & \xrightarrow{\pi_2} & M_2. \end{array}$$

Let $F = F_1 \oplus F_2$, so that $F$ is analytically projective again. Define

$$\pi \colon F \to M : (f_1, f_2) \mapsto \pi_1(f_1) + \psi_2(f_2).$$

It is an easy exercise in diagram chasing to show that $\pi$ is surjective. We claim that $\pi^{-1}(\mathfrak{p} M) = \mathfrak{p} F$, so that in particular $\ker \pi \subset \mathfrak{p} F$, proving that $\pi$ induces an analytically projective $\mathfrak{p}$-representation of $M$. So, to prove our claim, let $f = (f_1, f_2) \in F$ be such that $\pi(f) \in \mathfrak{p} M$. Hence

$$\begin{aligned} \theta(\pi(f)) &= \theta(\pi_1(f_1) + \psi_2(f_2)) \\ &= \theta \psi_2(f_2) \\ &= \pi_2(f_2) \in \mathfrak{p} M_2, \end{aligned}$$

which implies, since $\pi_2$ is an $\mathfrak{p}$-representation, that $f_2 \in \mathfrak{p} F_2$ and therefore that $\psi_2(f_2) \in \mathfrak{p} M$. But since $\pi(f) \in \mathfrak{p} M$, we have that $\pi_1(f_1) \in \mathfrak{p} M \cap M_1$. We now use that $M_1$ is $\mathfrak{p}$-closed in $M$, to obtain that $\pi_1(f_1) \in \mathfrak{p} M_1$ and using that also $\pi_1$ is an $\mathfrak{p}$-representation, we finally obtain that $f_1 \in \mathfrak{p} F_1$.

Now, the last statement is clear from the above construction. ■

## 3. Rigid Analytic Flatificators

**3.1. Definition.** Let $X$ be a rigid analytic variety. Given two locally closed analytic subvarieties $Z$ and $Z'$ of $X$, then we say that $Z \subset Z'$ (respectively, $Z = Z'$) at a point $x$ of $X$, if there exists an affinoid admissible open $U$ containing $x$, such that $Z \cap U \subset Z' \cap U$ (respectively, $Z \cap U = Z' \cap U$), or, more correctly, if the closed immersion $Z \cap U \hookrightarrow U$ factors through (respectively, coincides with) the closed immersion $Z' \cap U \hookrightarrow U$. Equivalently, if $\mathcal{I}$ and $\mathcal{I}'$ denote the coherent $\mathcal{O}_U$-ideals defining $Z$ and $Z'$, where $U$ is some admissible affinoid on which both closed analytic subvarieties are defined, then this is equivalent with $\mathcal{I}'_x \subset \mathcal{I}_x$ (respectively, $\mathcal{I}'_x = \mathcal{I}_x$).

A map of rigid analytic varieties $f \colon Y \to X$ is said to be flat in a point $y \in Y$, if the corresponding morphism of local rings $\mathcal{O}_{X, f(y)} \to \mathcal{O}_{Y, y}$ is flat and $f$ is called flat, if it is so in each point of $Y$. This is the standard way of viewing flatness as a property in the source (see for instance **(3.8)** below). However, we will also need to talk about flatness at a point of the target space $X$. Namely, we will say that $f$ is flat at $x \in X$, if, for any $y \in Y$ in the fibre of $x$ (i.e., such that $f(y) = x$), we have that the corresponding morphism of local rings $\mathcal{O}_{X, x} \to \mathcal{O}_{Y, y}$ is flat. The reader should check that this is equivalent with $\mathcal{O}_{X, x} \to (f_* \mathcal{O}_Y)_x$ being flat, where $f_* \mathcal{O}_Y$ denotes the $\mathcal{O}_X$-module obtained from $\mathcal{O}_Y$ by taking its direct image. In other words, for a map of affinoid varieties $f \colon \operatorname{Sp} B \to \operatorname{Sp} A = X$ to be flat at $x \in X$ is equivalent with the morphism $\mathcal{O}_{X, x} \to B \otimes_A \mathcal{O}_{X, x}$ to be flat.

Let $f\colon Y \to X$ be a map of rigid analytic varieties and let $x \in X$ with non-empty fibre $f^{-1}(x)$. We will say that $f$ has a *rigid analytic flatificator* at $x$, if there exists a locally closed analytic subvariety $Z$ of $X$ containing $x$, for which the following two conditions hold.

(1) The restriction map $f^{-1}(Z) \to Z$ is flat.
(2) For any locally closed analytic subvariety $V$ of $X$ containing $x$, such that the restriction $f^{-1}(V) \to V$ is flat at $x$, we have that $V \subset Z$ at $x$.

In other words, $Z$ is maximal (among all affinoid locally closed immersions) with respect to (1) at $x$. Even if a rigid analytic flatificator $Z$ exists, it need not be unique, but it is so at $x$: namely, any locally closed analytic subvariety $Z'$ containing $x$ such that (1) holds for it and such that $Z = Z'$ at $x$ will also be a flatificator of $f$ at $x$. Note also that we ask for more in (1) than only flatness in the point $x$. Therefore, if $f\colon Y = \operatorname{Sp} B \to X = \operatorname{Sp} A$ is a map of affinoid varieties, then our definition of flatificator translates to the ring morphism $\mathcal{O}_{X,x} \to (f_*\mathcal{O}_Y)_x$ having a flatificator $\mathfrak{f}$ plus the additional condition that $\mathcal{O}_X(U)/\mathfrak{f}\mathcal{O}_X(U) \to \mathcal{O}_Y(f^{-1}(U))/\mathfrak{f}\mathcal{O}_Y(f^{-1}(U))$ must be flat, for some affinoid subdomain $U$ of $X$ containing $x$ over which $\mathfrak{f}$ is defined. A rigid analytic flatificator at $x$ is then given by the closed analytic subvariety of $U$ defined by this ideal $\mathfrak{f}$.

We will say that $Z$ is a *universal rigid analytic flatificator* at $x$, if, for any map $\pi\colon \tilde{X} \to X$ of rigid analytic varieties and any $\tilde{x} \in \pi^{-1}(x)$, we have that $\pi^{-1}(Z)$ is a rigid analytic flatificator at $\tilde{x}$ of the base change map $Y \times_X \tilde{X} \to \tilde{X}$.

The theory of rigid analytic flatificators behaves only well for *quasi-compact* rigid analytic varieties, i.e., separated rigid analytic varieties admitting a finite admissible affinoid covering.

The following facts are easy to prove.

**3.2. Observation.** *If $Z$ is a rigid analytic flatificator at $x$ of a map $f\colon Y \to X$ of affinoid varieties, where $x \in X$ has non-empty fibre and if $W$ is a locally closed analytic subvariety of $X$ containing $x$, then $Z \cap W$ is a rigid analytic flatificator at $x$ of the restriction $f^{-1}(W) \to W$.*

*If $f\colon Y \to X$ is a map of rigid analytic varieties, $x \in X$ with non-empty fibre and $\{Y_i\}_{i<n}$ a finite admissible covering of $Y$, such that $f_i\colon Y_i \subset Y \to X$ has a flatificator $Z_i$ at $x$, for all $i < n$, then $Z = Z_0 \cap \cdots \cap Z_{n-1}$ is a flatificator of $f$ at $x$.*

As a corollary to the second observation, we see that in order to prove the existence of a universal rigid analytic flatificator when dealing with quasi-compact rigid analytic varieties, we can always reduce to the affinoid case. Furthermore, since any map $\pi\colon \tilde{X} \to X$ of affinoid varieties factors as a closed immersion $\tilde{X} \hookrightarrow X \times R^n$ and a projection $X \times R^n \twoheadrightarrow X$, for some $n$, we can, after applying the first observation, restrict ourselves in checking universality, to base changes over projections $\pi\colon \tilde{X} = X \times R^n \twoheadrightarrow X$.

**3.3. Theorem.** *Let $X = \operatorname{Sp} A$ be an affinoid variety and $x \in X$ an arbitrary point. Let $M$ be a finitely generated $A\langle T \rangle$-module. Then there exists an affinoid subdomain $\operatorname{Sp} C = U$ of $X$ containing $x$, such that $M \otimes_{A\langle T \rangle} C\langle T \rangle$ has a strongly founded, analytically projective and flat $\mathfrak{m}$-representation over $C$, where $\mathfrak{m}$ is the maximal ideal in $C$ corresponding to $x$.*

*Proof.* Since $M$ is a finite $A\langle T \rangle$-module, it is an analytic $A$-module. We will prove the theorem by induction on $n$, the number of variables $T = (T_1, \ldots, T_n)$. Let us first prove the case $n = 0$. Write $\mathfrak{m}$ as well for the maximal ideal in $A$ corresponding to $x$ and let

$$(1) \qquad 0 \to H' \longrightarrow A_{\mathfrak{m}}^d \xrightarrow{\pi'} M_{\mathfrak{m}} \to 0$$

be a minimal free resolution (and hence, in particular, a free $\mathfrak{m}A_{\mathfrak{m}}$-representation) of $M_{\mathfrak{m}}$. (The existence of such a representation for a finite module over a local ring is a well-known consequence of Nakayama's Lemma). In fact, there exists an $A$-module morphism $\pi\colon A^d \to M$ which localizes to $\pi'$. Let $Z$ denote the cokernel of $\pi$ and $H$ its kernel. By (1), there exists $s \in A \setminus \mathfrak{m}$, such that $sZ = 0$ and $sH \subset \mathfrak{m}A^d$. Let $\alpha = |s(x)| > 0$ and let $U = \operatorname{Sp} C \subset X$ be defined by $z \in U$, if and only if, $\alpha \le |s(z)|$. Note that $x \in U$. Since $s$ is a unit in $C$, we see that $Z \otimes_A C = 0$. Since $A \to C$ is flat, $H \otimes_A C$ is the kernel of the map $C^d \twoheadrightarrow M \otimes_A C$ obtained from $\pi$ by base change. Using once more that $s$ is a unit, we conclude that $H \otimes_A C \subset \mathfrak{m}C^d$, in other words, that

$$0 \to H \otimes_A C \longrightarrow C^d \longrightarrow M \otimes_A C \to 0$$

is a free $\mathfrak{m}$-representation of the latter module.

So, assume $n > 0$ and the theorem proven for $n - 1$ and for all affinoid algebras. Let $\mathfrak{P} = \mathfrak{m}A\langle T\rangle$, so that $\mathfrak{P}$ is a prime ideal of $A\langle T\rangle$. Using now **(1.6)** and **(2.5)** on the prime ideal $\mathfrak{P}$ of $A\langle T\rangle$, we can reduce to the following two cases that either $M$ is cyclic or $M_\mathfrak{P} = 0$. Just use that any $\mathfrak{P}$-series given as in **(1.6)** remains of the same type after a flat base change $A \to C$, that is, after restricting it to an affinoid subdomain. Note also that any surjection appearing here is automatically strict, due to **(2.2)**.

*Case 1.* Suppose that $M$ is cyclic, hence, as an $A\langle T\rangle$-module, isomorphic with $A\langle T\rangle/I$, where $I$ is an ideal in $A\langle T\rangle$. If $I \subset \mathfrak{P}$, then the short exact sequence

$$0 \to I \longrightarrow A\langle T\rangle \longrightarrow M \to 0$$

is an $\mathfrak{m}$-representation of the desired type by **(1.2)** and **(2.4)**. If $I \not\subset \mathfrak{P}$, then $M_\mathfrak{P} = 0$ and we are in case two.

*Case 2.* Suppose that $M_\mathfrak{P} = 0$. So there exists an $s \notin \mathfrak{P}$, such that $sM = 0$. Hence $M$ is an $A\langle T\rangle/(s)$-module. Write $s$ as $\sum a_i T^i$, with $a_i \in A$. Since $s \notin \mathfrak{P} = \mathfrak{m}A\langle T\rangle$, there is some coefficient $a_i$ which does not vanish at $x$. In particular, the maximum $\alpha$ of all $|a_i(x)|$ is non-zero. Let $\nu$ be the lexicographically largest index $i$ for which $|a_\nu(x)| = \alpha$. Since $|a_i|$ goes to zero when $i$ goes to infinity, we can find an $\varepsilon > 0$, such that, for all $i > \nu$, we have that $|a_i(x)| \leq \alpha - \varepsilon$. Let $U = \operatorname{Sp} C$ be the affinoid subdomain of $X$ given by $z \in U$, if and only if,

(2)
$$|a_\nu(z)| \geq \alpha$$
$$|a_i(z)| \leq \alpha - \varepsilon \qquad \text{for } i > \nu.$$

Then $x \in U$ and $a_\nu$ is a unit in $C$ (since it vanishes nowhere). Moreover, for any $z \in U$, we have that $|a_i(z)/a_\nu(z)| < 1$, whenever $i > \nu$. Therefore, we may assume after applying a (Weierstrass) automorphism that $s$ as an element of $C\langle T\rangle$ is regular in $T_n$. By the Weierstrass Preparation Theorem, we then have that $C\langle T\rangle/(s)$ is a finitely generated $C\langle T'\rangle$-module, where $T' = (T_1, \ldots, T_{n-1})$. So, $M \otimes_{A\langle T\rangle} C\langle T\rangle$ is a finitely generated $C\langle T'\rangle$-module and by induction we have that $M \otimes_{A\langle T\rangle} C\langle T\rangle$ has an $\mathfrak{m}$-representation as desired. ∎

*Remark.* Note that by [**3**, 3.7.3. Proposition 6], we have that

$$M \otimes_{A\langle T\rangle} C\langle T\rangle = M \widehat{\otimes}_{A\langle T\rangle} C\langle T\rangle.$$

Note also that nowhere we made use of the fact that $\mathfrak{m}$ was a maximal ideal, but only that it was prime. This will be all what is needed in the Berkovich equivalent of this proposition as discussed in the Appendix.

**3.4. Lemma.** *Let $A$ be an affinoid algebra and $\mathfrak{a}$ an ideal of $A$. Let $g \in A\langle S\rangle$, with $S$ a single variable. If there exists a sequence $(r_m)_{m \in \mathbb{N}}$ of non-zero elements of $R$ converging to zero such that $g(r_m) \in \mathfrak{a}$, for all $m$, then $g \in \mathfrak{a}A\langle S\rangle$.*

*Proof.* By considering $A/\mathfrak{a}$ instead of $A$, we may reduce to the case $\mathfrak{a} = 0$. Endow $A$ with some residue norm. Suppose $g \neq 0$, but all $g(r_m) = 0$. Write

$$g(S) = \sum_{d=i}^{\infty} a_i S^i,$$

where $a_i \in A$ with $a_d \neq 0$. Let

$$\alpha = \max_{i > d} |a_i|.$$

There is nothing to prove if $\alpha = 0$ (since a polynomial cannot have infinitely many roots unless it is the zero polynomial), so assume $\alpha > 0$. Choose $m \gg 0$, such that $|r_m| < |a_d|/\alpha$. Since $g(r_m) = 0$ and $r_m \neq 0$, it follows that $|a_d| = \left|\sum_{i>d} a_i r_m^{i-d}\right|$. On the other hand, for $i > d$, we have that

$$\left|a_i r_m^{i-d}\right| \leq |a_i| \cdot |r_m| \leq \alpha |r_m| < |a_d|,$$

contradiction. ∎

**3.5. Lemma.** *Let $A$ be an affinoid algebra and $\mathfrak{p}$ a prime ideal of $A$. Let $f \in A\langle S\rangle$, with $S$ a single variable. If there exists a sequence $(r_m)_{m\in\mathbb{N}}$ of non-zero elements of $R$ converging to zero such that $f(r_m) = 0$ in $A_\mathfrak{p}$, then $f = 0$ in $A\langle S\rangle_\mathfrak{p} = A\langle S\rangle \otimes_A A_\mathfrak{p}$ (i.e., the localization at $\mathfrak{p}$ of the $A$-module $A\langle S\rangle$, which in general is smaller than the localization of $A\langle S\rangle$ at its prime ideal $\mathfrak{p}A\langle S\rangle$).*

*Proof.* By induction, we will show that $f \in \mathfrak{p}^n A\langle S\rangle_\mathfrak{p}$, for all $n \in \mathbb{N}$, from which the lemma then follows immediately by Krull's Intersection Theorem, since $\mathfrak{p}$ lies in the radical of $A\langle S\rangle_\mathfrak{p}$.

Let $Q$ denote the ring $A\langle S\rangle_\mathfrak{p}$. Suppose that $f \in \mathfrak{p}^n Q$, so that we can write

$$(1) \qquad f(S) = \sum_{i<t} \mu_i f_i(S),$$

in $Q$, where the $\mu_i$ form a minimal set of generators of $\mathfrak{p}^n$ in $A_\mathfrak{p}$ and $f_i(S) \in A\langle S\rangle$. By assumption, we have, for each $m$, that

$$0 = f(r_m) = \sum_{i<t} \mu_i f_i(r_m)$$

in $A_\mathfrak{p}$. Since the $\mu_i$ minimally generate $\mathfrak{p}^n A_\mathfrak{p}$, we must have that each $f_i(r_m) \in \mathfrak{p}A_\mathfrak{p}$ and whence $f_i(r_m) \in \mathfrak{p}$. Our claim now follows after applying **(3.4)** with $g = f_i$ and $\mathfrak{a} = \mathfrak{p}$. ∎

**3.6. Proposition.** *Let $X$ be an affinoid variety and $x \in X$. Let $V \subset X \times R$ be a closed analytic subvariety, containing $(x,0)$. Let $(r_m)_m$ be a sequence of non-zero elements of $R$ tending to zero and let $\sigma_m$ be the automorphism (translation) of $X \times R$ given by $(y,a) \mapsto (y, a + r_m)$. If, for all $m$, we have that $\sigma_m(V) \subset V$ at $(x,0)$, then there exists a locally closed analytic subvariety $W$ of $X$ containing $x$, such that $V = W \times R$ at $(x,0)$.*

*Proof.* Let $X = \operatorname{Sp} A$ and let $\mathfrak{m}$ be the maximal ideal of $A$ associated to $x$. Hence $A\langle S\rangle$ is the associated affinoid algebra of $X \times R$, where $S$ is a single variable. Let $\mathfrak{G}$ be the ideal of $A\langle S\rangle$ defining $V$ and put $B = A\langle S\rangle / \mathfrak{G}$. Let $\mathfrak{M} = (\mathfrak{m}, S)$ be the maximal ideal of $A\langle S\rangle$ corresponding to $(x,0)$. Let $(G_1, \ldots, G_s) = \mathfrak{G}$, with $G_i \in A\langle S\rangle$. Introduce a new variable $T$ and let $P_i(S,T) = G_i(S-T)$. In order to obtain a proof, we will exploit the fact that $P_i$ has a Taylor expansion with respect to the variable $T$, and whence with coefficients functions of $S$ (corresponding in characteristic zero to the derivatives of $G_i$). That is to say, we can write

$$(1) \qquad P_i(S,T) = G_i(S-T) = \sum_j h_{ij}(S) T^j,$$

with $h_{ij} \in A\langle S\rangle$. By assumption, we have, for each $m$, that

$$\sigma_m(G_i) = G_i(S - r_m) = P_i(S, r_m) \in \mathfrak{G} A\langle S\rangle_\mathfrak{M}.$$

In other words, we see that $P_i(S, r_m) = 0$ when considered as an element of $B_{\mathfrak{M}B}$. Applying **(3.5)**, we obtain that all $P_i \in B\langle T\rangle_{\mathfrak{M}B}$. Therefore, there exists some $Q(S) \in A\langle S\rangle$ not contained in $\mathfrak{M}$, such that $Q(S) G_i(S - T) \in \mathfrak{G} A\langle S, T\rangle$, for $i = 1, \ldots, s$. Hence, by considering the coefficients of the $T^j$ using (1), we have that $Q(S) h_{ij}(S) \in \mathfrak{G}$, for all the relevant $i$ and $j$. Since $Q(x,0) \neq 0$ by assumption, we can find an admissible affinoid $U' = \operatorname{Sp}(C')$ of $X \times R$ containing $(x,0)$, such that all $h_{ij} \in \mathfrak{G}C'$. By shrinking $U'$ further (we are only seeking to prove something locally around $(x,0)$), we can assume that $U'$ is of the form $U \times \mathbb{B}$, where $U = \operatorname{Sp} C$ is an affinoid subdomain of $X$ containing $x$ and $\mathbb{B}$ a disk of small radius around 0. The latter is isomorphic to $R$ by a map of the form $S \mapsto \pi S$. Since the image of the sequence $(r_m)_m$ remains a zero sequence under such a map, we may assume without loss of generality that $U' = U \times R$. Hence $h_{ij} \in \mathfrak{G}C\langle S\rangle$, for all $i$ and $j$. Substituting $S - r_m$ for $T$ in (1), we have, using [**3**, 5.2.7. Corollary 2] that

$$(2) \qquad G_i(r_m) = \sum_j h_{ij}(S)(S - r_m)^j \in \mathfrak{G}C\langle S\rangle.$$

Let $\mathfrak{g} = \mathfrak{G}C\langle S\rangle \cap C$, so that actually $G_i(r_m) \in \mathfrak{g}$, for all $i$ and $m$. Applying **(3.4)**, we conclude that $G_i \in \mathfrak{g}C\langle S\rangle$, proving that $\mathfrak{g}C\langle S\rangle = \mathfrak{G}C\langle S\rangle$, from which the proposition now follows readily. ∎

**3.7. Theorem.** *Let $f\colon Y \to X$ be a map of quasi-compact rigid analytic varieties and let $x$ be an arbitrary point of $X$ with non-empty fibre $f^{-1}(x)$. Then there exists an admissible affinoid $U$ of $X$ containing $x$, such that the restriction $f^{-1}(U) \to U$ has a universal rigid analytic flatificator at $x$.*

*Proof.* Since this is a local question, we may already assume, by our remarks in **(3.2)**, that $X = \operatorname{Sp} A$ and $Y = \operatorname{Sp} B$ are affinoid. From the induced algebra morphism $A \to B$, we can construct a surjective algebra morphism $A\langle T\rangle \to B$, turning $B$ into a finite $A\langle T\rangle$-module, for some variables $T = (T_1, \ldots, T_n)$. Applying **(3.3)** to this situation, gives the existence of an affinoid subdomain $U = \operatorname{Sp} C$ around $x$, such that

$$D = \mathcal{O}_Y(f^{-1}(U)) = B\widehat{\otimes}_A C = B \otimes_{A\langle T\rangle} C\langle T\rangle$$

has a strongly founded and flat $\mathfrak{m}$-representation

(1) $$0 \to H \longrightarrow F \longrightarrow D \to 0,$$

where $\mathfrak{m}$ is the maximal ideal of $C$ corresponding to $x$. Taking stalks at the point $x$, we obtain a $\mathfrak{m}\mathcal{O}_{X,x}$-representation

$$0 \to H_x \longrightarrow F_x \longrightarrow D_x \to 0,$$

where we have written $D_x$ for $(f_*\mathcal{O}_Y)_x$. From **(1.1)** it then follows that $C_x = \mathcal{O}_{X,x} \to D_x$ has a flatificator $\mathfrak{f}$. By shrinking $U$ further if necessary, we can assume that $\mathfrak{f}$ is actually an ideal of $C$ and that $H \subset \mathfrak{f}F$. Therefore, it follows immediately from (1) that $D/\mathfrak{f}D$ is flat. As observed in **(3.1)**, the locally closed analytic subvariety given by $\mathfrak{f}$ is then a rigid analytic flatificator of $f$ at $x$.

Hence the theorem will be proven, if we can show that any rigid analytic flatificator is universal. To this end, let $f\colon Y \to X$ be a map of affinoid varieties with rigid analytic flatificator $Z$ and let $\pi\colon X \times R^n \to X$ be the canonical projection map. We need to show that, for each $r \in R^n$, we have that $Z \times R^n$ is a rigid analytic flatificator of $f \times 1\colon Y \times R^n \to X \times R^n$ in $(x, r)$. By an inductive argument we may reduce to the case $n = 1$ and after a translation, we can assume that $r = 0$.

First of all, since $f^{-1}(Z) \to Z$ is flat, also its base change

(2) $$f^{-1}(Z) \times R \to Z \times R \quad \text{is flat.}$$

By above we know that there exists an affinoid subdomain $U$ of $X \times R$ and a closed analytic subvariety $V$ of $U$, both containing $(x, 0)$, such that $V$ is a rigid analytic flatificator of $(f \times 1)^{-1}(U) \to U$ at $(x, 0)$. Since a flatificator is only locally defined, it follows from (2) that we only need to show that $Z = V$ at $(x, 0)$. Hence without loss of generality we may assume that $V$ is a closed analytic subvariety of $X \times R$ (i.e., is globally defined). From (2) and the definition of $V$, it follows that $Z \times R \subset V$ at $(x, 0)$ and we seek to prove the opposite inclusion.

From the latter inclusion it also follows that we can find a sequence $(r_m)_m$ of non-zero elements in $R$ converging to zero, such that the associated translation $\sigma_m\colon (y, a) \mapsto (y, a + r_m)$ has the property that $(x, 0) \in \sigma_m(V)$. Since $\sigma_m$ is an automorphism and since $(f \times 1)^{-1}(V) \to V$ by definition is flat, we have that the map

$$(f \times 1)^{-1}(\sigma_m(V)) \to \sigma_m(V)$$

is also flat. By definition of rigid analytic flatificator, this implies that $\sigma_m(V) \hookrightarrow V$ at $(x, 0)$.

We can now invoke **(3.6)**, to find a locally closed analytic subvariety $W$ of $X$ containing $x$, such that $V = W \times R$ at $(x, 0)$. In particular, the map

(3) $$f^{-1}(W) \times R \to W \times R$$

is then flat at $(x, 0)$. Since the map $W \times R \to W$ is faithfully flat, we obtain from (3) that then also the map $f^{-1}(W) \to W$ must be flat at $x$. Hence, by definition of $Z$, we have that $W \subset Z$ at $x$, proving that $V = W \times R \subset Z \times R$ at $(x, 0)$. ∎

Next we list some corollaries to the existence of (universal) rigid analytic flatificators. The first is the rigid analytic analogue of [**4**, $O_{III}$ 11.1.1].

**3.8. Theorem.** *Let $f\colon Y \to X$ be a map of rigid analytic varieties and let $y \in Y$. If $f$ is flat at $y$, then there exists an admissible affinoid open $U$ of $Y$ containing $y$, such that the restriction $f|_U \colon U \to X$ is flat.*

*Proof.* Since this is a local question, we may already assume that $X = \operatorname{Sp} A$ and $Y = \operatorname{Sp} B$ are affinoid. Let $\mathfrak{n}$ be the maximal ideal of $B$ corresponding to the point $y$ and let $\mathfrak{m} = \mathfrak{n} \cap A$, which then corresponds to the point $x$. View the field $K$ as an $A$-module via the isomorphism $K \cong A/\mathfrak{m}$. Let $\mathcal{F}$ be the coherent $\mathcal{O}_Y$-module defined by the finite $B$-module $\operatorname{Tor}_1^A(B,K)$. (Note that the latter module is indeed finite, since it is the kernel of the morphism $\mathfrak{m} \otimes_A B \to B$.) Since $B_\mathfrak{n}$ is flat over $A$, we have that $\mathcal{F}_y = 0$. Hence there exists an admissible affinoid $U = \operatorname{Sp} C$ of $Y$ containing $y$, such that $\mathcal{F}(U) = 0$. In other words, we obtained that $\operatorname{Tor}_1^A(C,K) = 0$.

Applying **(3.3)** to the morphism $A \to C$, we can find an admissible open $W = \operatorname{Sp} A'$ of $X$ containing $x$, such that $C' = C \widehat{\otimes}_A A' = \mathcal{O}_Y(f^{-1}(W))$ has a founded and flat $\mathfrak{m}A'$-representation

$$(1) \qquad 0 \to H \longrightarrow F \longrightarrow C' \to 0,$$

where $H \subset \mathfrak{m}F$ and $F$ is a flat and founded $A'$-module. Since $\operatorname{Tor}_1^A(C,K) = 0$ and $A \to A'$ is flat, we also have that $\operatorname{Tor}_1^{A'}(C',K) = 0$.[1] After tensoring (1) with $K$ over $A'$, we then conclude that $H \otimes_{A'} K = 0$, i.e., that $H = \mathfrak{m}H$. Hence, by induction, we see that $H \subset \mathfrak{m}^n F$, for all $n$. Let

$$F^{(\infty)} = \bigcap_n \mathfrak{m}^n F,$$

so that $H \subset F^{(\infty)}$. Since $F$ is $A'$-founded, we must have that

$$F^{(\infty)} \subset \left(\bigcap_n \mathfrak{m}^n A'\right) F.$$

But the stalk at $x$ of the latter module is zero, since $x \in W$, so that $\mathcal{O}_{W,x} = \mathcal{O}_{X,x}$ is $\mathfrak{m}$-adically separated. Hence also $H_x = 0$ so that after shrinking $W$ even further we may already assume that $H = 0$. Therefore we proved that $C' = F$ is $A'$-flat whence $A$-flat, proving that $f^{-1}(W) = \operatorname{Sp} C'$ is the required admissible affinoid. ∎

**3.9. Theorem.** *Let $f\colon Y \to X$ be a map of quasi-compact rigid analytic varieties. Suppose that there exist maps of rigid analytic varieties $\pi_i \colon X_i \to X$ with the following two properties.*
  (1) *For each $i$, the base change $Y \times_X X_i \to X_i$ is flat.*
  (2) *The $\{\pi_i\}_i$ form a surjective family, i.e.,*

$$\bigcup_i \pi_i(X_i) = |X|.$$

*If $X$ is reduced, then $f$ is flat.*

*Proof.* Let $x \in X$ and let $U$ be an admissible affinoid of $X$ containing $x$, such that the restriction $f^{-1}(U) \to U$ has a universal rigid analytic flatificator $Z$ at $x$, the existence of which is assured by **(3.7)**. If $U_i$ denotes $\pi_i^{-1}(U)$, then the base change $f^{-1}(U) \times_U U_i \to U_i$ is also flat, by (1), so that by universality of $Z$, we must have that $\pi_i^{-1}(Z) = U_i$ at any point of the fibre of $x$. Hence, by (2), we therefore have that $|Z| = |X|$ at $x$. Since $X$ is reduced, we even have that $X = Z$ at $x$ and whence that $f$ is flat at $x$. Since this holds for any point $x \in X$, we proved that $f$ is flat. ∎

---

[1]To see this, let $P_\bullet$ be a finitely generated projective resolution of $K$ as an $A$-module. Since $A \to A'$ is flat, we obtain that $P_\bullet \otimes_A A'$ is a projective resolution of $K$ as an $A'$-module. Hence

$$(a) \qquad \operatorname{Tor}_i^{A'}(C',K) \cong H_i((P_\bullet \otimes_A A') \otimes_{A'} C').$$

Considering $C$ as some quotient of $A\langle T \rangle$, we see that $C' = C \otimes_{A\langle T \rangle} A'\langle T \rangle$. Therefore, the right hand side of (a) equals

$$H_i((P_\bullet \otimes_A C) \otimes_{A\langle T \rangle} A'\langle T \rangle)$$

which by flatness equals

$$H_i(P_\bullet \otimes_A C) \otimes_{A\langle T \rangle} A'\langle T \rangle \cong \operatorname{Tor}_i^A(C,K) \otimes_{A\langle T \rangle} A'\langle T \rangle.$$

4. Analytic Complements

**4.1. Definition.** Let $A$ be a ring and $\mathfrak{a}$ an ideal of $A$. We introduce the following ideal

$$H_A(\mathfrak{a}) = \sum_{n=1}^{\infty} \operatorname{Ann}_A(\mathfrak{a}^n).$$

Note that if $A$ is Noetherian, then there exists some $n$, depending on the ideal $\mathfrak{a}$, such that $H_A(\mathfrak{a}) = \operatorname{Ann}_A(\mathfrak{a}^n)$. If $A$ is reduced, then $H_A(\mathfrak{a}) = \operatorname{Ann}_A(\mathfrak{a})$. Let $X$ be a rigid analytic variety and let $\mathcal{I}$ be a coherent $\mathcal{O}_X$-ideal. We define similarly

$$\mathcal{H}_X(\mathcal{I}) = \sum_{n=1}^{\infty} \mathcal{A}nn_{\mathcal{O}_X}(\mathcal{I}^n).$$

As observed in the proof of [**18**, Proposition 2.2.1], this is again a coherent $\mathcal{O}_X$-ideal and, moreover, if $U = \operatorname{Sp} A$ is an admissible affinoid of $X$ and $\mathfrak{a} = \mathcal{I}(U)$, then $\mathcal{H}_X(\mathcal{I})(U) = H_A(\mathfrak{a})$.

Let $X$ be a rigid analytic variety and let $Z$ be a closed analytic subvariety, given by the coherent ideal $\mathcal{I}$. We define the *analytic complement* of $Z$ as the closed analytic subvariety defined by $\mathcal{H}_X(\mathcal{I})$ and denote it $\complement_X(Z)$.

**4.2. Lemma.** *Let $X$ be a rigid analytic variety and let $Z$ be a closed analytic subvariety of $X$. Then, as point sets,*

(1) $$\operatorname{cl}(X \setminus Z) = |\complement_X(Z)|.$$

*Remark.* The *(Zariski-)closure* $\operatorname{cl}(U)$ of a subset $U$ of $X$ is defined as the smallest closed analytic subset containing $U$.

*Proof.* Let $\mathcal{I}$ be the coherent $\mathcal{O}_X$-ideal defining $Z$ and let $\mathcal{H} = \mathcal{H}_X(\mathcal{I})$. We claim that $\mathcal{I} \cdot \mathcal{H}$ is nilpotent and hence that $|Z| \cup |\complement_X(Z)| = |X|$, i.e., the $\subset$-inclusion in (1) holds. It is enough to verify our claim in each stalk. Hence let $x \in X$, so that $\mathcal{H}_x = \operatorname{Ann}_{\mathcal{O}_{X,x}}(\mathcal{I}_x^n)$, for some $n$. But then clearly $(\mathcal{I}_x \mathcal{H}_x)^n = 0$.

The opposite inclusion will follow, if we can show that any closed analytic subset $V$ containing $X \setminus Z$, contains also $|\complement_X(Z)|$. Therefore, using that $|Z| \cup V = |X|$ we obtain that $\mathcal{I} \cdot \mathcal{J}$ is nilpotent, where $\mathcal{J}$ denotes the (radical) coherent $\mathcal{O}_X$-ideal defining $V$. Hence, at each point $x \in X$ there is some $n$, such that $\mathcal{J}_x^n \subset \mathcal{H}_x$. But then $\mathcal{J}$ is contained in the radical of $\mathcal{H}$, which is equivalent with $|\complement_X(Z)| \subset V$, as we wanted to show. ∎

**4.3. Lemma.** *Let $X = \operatorname{Sp} A$ be an affinoid variety and let $\mathfrak{a}, \mathfrak{f}_1, \mathfrak{f}_2$ be ideals of $A$. Let $\pi_i \colon \tilde{Y}_i \to Y$ be the blowing up of $Y$ with respect to $\mathfrak{f}_i$, for $i = 1, 2$. Let $V$ be the closed analytic subvariety of $X$ defined by $\mathfrak{a}$. If*

$$\mathfrak{f}_1 + H_A(\mathfrak{a}) = \mathfrak{f}_2 + H_A(\mathfrak{a}),$$

*then*

$$\complement_{\tilde{Y}_1}(\pi_1^{-1}(V)) \cong \complement_{\tilde{Y}_2}(\pi_2^{-1}(V)).$$

*Moreover, if $\mathfrak{f}_i \subset \operatorname{rad}(\mathfrak{a})$, then $\complement_{\tilde{Y}_i}(\pi_i^{-1}(V)) = \tilde{Y}_i$.*

*Proof.* Let us denote $\tilde{W}_i = \complement_{\tilde{Y}_i}(\pi_i^{-1}(V))$, which, by definition is the closed analytic subvariety of $\tilde{Y}_i$ defined by $\mathcal{H}_i = \mathcal{H}_{\tilde{Y}_i}(\mathfrak{a}\mathcal{O}_{\tilde{Y}_i})$. We claim that $\mathfrak{f}_2 \mathcal{O}_{\tilde{W}_1}$ is invertible. To see this, let $U_1 = \operatorname{Sp} B_1$ be an admissible affinoid of $\tilde{Y}_1$. Hence, if we put $C_1 = B_1/H_{B_1}(\mathfrak{a} B_1)$, then $C_1 = \mathcal{O}_{\tilde{W}_1}(U_1 \cap \tilde{W}_1)$. Since $H_A(\mathfrak{a}) \, C_1 = 0$, we see that $\mathfrak{f}_1 C_1 = \mathfrak{f}_2 C_1$. Since $\mathfrak{f}_1 B_1$ is invertible, by definition of a blowing up, we have that $\mathfrak{f}_1 C_1$ is locally principal. So, in order to prove our claim, we have to verify that $\operatorname{Ann}_{C_1}(\mathfrak{f}_1 C_1) = 0$. Hence let $b \in B_1$, such that $b\mathfrak{f}_1 C_1 = 0$. This means that $b\mathfrak{f}_1 \subset H_{B_1}(\mathfrak{a} B_1)$ and hence that $b\mathfrak{f}_1 \mathfrak{a}^n B_1 = 0$, for some $n$. But $\mathfrak{f}_1$ is invertible in $B_1$, so that $b\mathfrak{a}^n B_1 = 0$, from which it follows that $b = 0$ as an element of $C_1$.

Therefore, we proved that $\mathfrak{f}_2 \mathcal{O}_{\tilde{W}_1}$ is invertible and hence, there exists a unique map of rigid analytic varieties $i_1 \colon \tilde{W}_1 \to \tilde{Y}_2$ making the following diagram commute

(1)
$$\begin{array}{ccc} \tilde{W}_1 & \xrightarrow{\pi_1|_{\tilde{W}_1}} & Y \\ {\scriptstyle i_1}\downarrow & & \| \\ \tilde{Y}_2 & \xrightarrow{\pi_2} & Y. \end{array}$$

We claim that $i_1$ factors through $\tilde{W}_2$, i.e., that $i_1(\tilde{W}_1) \subset \tilde{W}_2$. In order to prove this, we must show that $\mathcal{H}_2 \mathcal{O}_{\tilde{W}_1} = 0$. Since this is a local question, we may take an admissible affinoid $U_2 = \operatorname{Sp} B_2$, such that $i_1(U_1 \cap \tilde{W}_1) \subset U_2$ and we need to show that $H_{B_2}(\mathfrak{a}B_2)C_1 = 0$, under the induced morphism $B_2 \to C_1$ of affinoid algebras. Hence let $b_2 \in H_{B_2}(\mathfrak{a}B_2)$, so that there exists some $n$, for which $b_2\mathfrak{a}^n B_2 = 0$. Hence also $b_2\mathfrak{a}^n C_1 = 0$, which implies that, for some $m$, we have that $b_2\mathfrak{a}^{m+n} B_1 = 0$. But the latter means that $b_2 \in H_{B_1}(\mathfrak{a}B_1)$, i.e., that $b_2 = 0$ as an element of $C_1$, as we needed to show.

By symmetry, we find a unique map $i_2 \colon \tilde{W}_2 \to \tilde{Y}_1$, making a similar diagram as (1) commute in which the indices 1 and 2 are interchanged, and such that $i_2(\tilde{W}_2) \subset \tilde{W}_1$. Uniqueness of both maps then implies that they must be each others inverse, providing the required isomorphism.

To prove the last assertion, we need to show that $\mathcal{H}_i = 0$ under the additional assumption that $\mathfrak{f}_i \subset \operatorname{rad}(\mathfrak{a})$. Let $m$ be such that $\mathfrak{f}_i^m \subset \mathfrak{a}$. By symmetry, we only have to deal with the case $i = 1$. We can calculate this locally, so that with the notation of above, we need to show that $H_{B_1}(\mathfrak{a}B_1) = 0$. Hence let $s$ lie in the latter ideal, so that $s\mathfrak{a}^n B_1 = 0$, for some $n$. But this implies by our assumption that also $s\mathfrak{f}_1^{mn} B_1 = 0$ and hence that $s = 0$ in $B_1$, since $\mathfrak{f}_1 B_1$ is invertible. ∎

**4.4. Corollary.** *Let $X = \operatorname{Sp} A$ be an affinoid variety and let $V$ and $Z_1$ be two closed analytic subvarieties of $X$, such that $|V| \subset |Z|_1$. Let $Z_2$ be another closed analytic subvariety of $X$, such that we have closed immersions $Z_1 \cap \complement_X(V) \subset Z_2 \subset Z_1$. Let $\pi_i \colon \tilde{X}_i \to X$ be the blowing up of $X$ with centre $Z_i$, for $i = 1, 2$. Then*
$$\pi_1^{-1}(\complement_X(V)) \cong \pi_2^{-1}(\complement_X(V)).$$

*Proof.* Let $\mathfrak{a}$ be the ideal defining $V$, so that $\complement_X(V)$ is defined by the ideal $H_A(\mathfrak{a})$. From our assumptions it follows that $Z_1 \cap \complement_X(V) = Z_2 \cap \complement_X(V)$, so that $\mathfrak{f}_1 + H_A(\mathfrak{a}) = \mathfrak{f}_2 + H_A(\mathfrak{a})$, where $\mathfrak{f}_i$ is the ideal defining $Z_i$, for $i = 1, 2$. Moreover, we have that $\mathfrak{f}_1 \subset \operatorname{rad}(\mathfrak{a})$. Hence from **(4.3)** we deduce the existence of an isomorphism
$$\tilde{X}_1 \cong \complement_{\tilde{X}_2}(\pi_2^{-1}(V)).$$
As is easily seen, we have that $\pi_2^{-1}(\complement_X(V))$ lies in the latter space, from which the required isomorphism follows readily. ∎

**4.5. Proposition.** *Let $X = \operatorname{Sp} A$ be an affinoid variety and let $Z$ be a closed analytic subvariety of $X$. Then the following are equivalent.*
   (i) *$Z$ is nowhere dense in $X$.*
   (ii) *$|\complement_X(Z)| = |X|$.*
   (iii) *The pair $(X, Z)$ satisfies the minimality condition.*
   (iv) *The kernel of the restriction map $A \to \mathcal{O}_X(X \setminus Z)$ is nilpotent.*
   (v) *$|Z|$ does not contain an irreducible component of $X$.*

*Remark.* Here we say that a subset $\Gamma$ is *nowhere dense* in $X$, if $X \setminus \Gamma$ is dense in $X$ with respect to the Zariski topology on $X$, i.e., each non-empty Zariski-open $U$ of $X$ intersects $X \setminus \Gamma$ non trivially. Clearly any subset of a nowhere dense subset is nowhere dense. Recall also that the pair $(X, Z)$ satisfies the minimality condition, if for any closed analytic subset $V$ of $X$ for which $X \setminus Z = V \setminus Z$, we must have that $V = |X|$. Strictly speaking, this is the set-theoretical version of the concept; for the more general notion, we require an equality of spaces $V = X$ instead. See [**6**, Proposition 3.8] for an analogous proposition using this more general notion.

*Proof.* Let $\mathfrak{a}$ be the ideal defining $Z$. By definition, $\complement_X(Z)$ is then defined by the ideal $H_A(\mathfrak{a})$.
   (i) $\implies$ (ii). Immediate from **(4.2)**.

$(ii) \implies (iii)$. Let $V$ be a closed analytic subset of $X$ such that $V \setminus Z = X \setminus Z$. Taking the Zariski closure, together with **(4.2)**, we obtain that
$$|X| = |\mathcal{C}_X(Z)| = \text{cl}(X \setminus Z) = \text{cl}(V \setminus Z) \subset V.$$

$(iii) \implies (iv)$. Let $f \in A$ be in the kernel of the restriction map $A \to \mathcal{O}_X(X \setminus Z)$. Let $V = V(f)$, then one verifies easily that $X \setminus Z = V \setminus Z$. By assumption, this implies that $V = |X|$ and hence that $f$ is nilpotent.

$(iv) \implies (i)$. Suppose $Z$ is not nowhere dense, so that there exists a proper closed analytic subset $V = V(\mathfrak{n})$ of $X$, such that $X \setminus V \subset Z$. This implies that $X \setminus Z = V \setminus Z$. But the restriction map $\rho \colon A \to \mathcal{O}_X(X \setminus Z)$ factors then through $\mathcal{O}_V(V) = A/\mathfrak{n}$, so that $\mathfrak{n}$ (which is not nilpotent) is contained in the kernel of $\rho$, contradiction.

$(iii) \implies (v)$. Suppose that $\Sigma$ is an irreducible component of $X$ such that $\Sigma \subset |Z|$. Let $X_0$ denote the union of all irreducible components of $X$ other than $\Sigma$. Hence we have that $X_0 \subsetneq |X|$ and $|X| = \Sigma \cup X_0$. But $|Z|$ contains $\Sigma$, so that $|X| = |Z| \cup X_0$. Hence by our assumption, we must have $|X| = X_0$, contradiction.

$(v) \implies (ii)$. Suppose that $(ii)$ does not hold. This means that $H_A(\mathfrak{a})$ is not nilpotent. Hence there exists a minimal prime $\mathfrak{p}$ of $A$, such that $H_A(\mathfrak{a}) \not\subset \mathfrak{p}$. Let $t \in H_A(\mathfrak{a})$ with $t \notin \mathfrak{p}$. Hence, for some $n$, we have that $t\mathfrak{a}^n = 0$, from which it follows that $\mathfrak{a}^n \subset \mathfrak{p}$ and hence $\mathfrak{a} \subset \mathfrak{p}$. This means that $V(\mathfrak{p}) \subset |Z|$, contradicting our assumption, since $V(\mathfrak{p})$ is an irreducible component of $X$. ∎

**4.6. Corollary.** *Let $\pi \colon \tilde{X} \to X$ be a blowing up map of a rigid analytic variety $X$ with a nowhere dense centre. If $X$ is reduced, then also $\tilde{X}$ is. Moreover, $\pi$ is surjective.*

*Proof.* Let $Z$ denote the centre of the blowing up $\pi$. Let $\tilde{\mathcal{G}} = \text{nil}(\mathcal{O}_{\tilde{X}})$. By [**18**, Theorem 3.2.1], the map $\pi$ is proper. Hence by [**3**, 9.6.3. Theorem 1], $\pi_*(\tilde{\mathcal{G}})$ is coherent and by the left exactness of $\pi_*$ it is an $\mathcal{O}_X$-ideal. We claim that $\pi_*(\tilde{\mathcal{G}})(X \setminus Z) = 0$, so that by **(4.5.iv)**, using that $X$ is reduced, we have that $\pi_*(\tilde{\mathcal{G}}) = 0$, from which it follows immediately that then also $\tilde{\mathcal{G}} = 0$, as we wanted to show. To prove our claim, observe that
$$\pi_*(\tilde{\mathcal{G}})(X \setminus Z) = \tilde{\mathcal{G}}(\tilde{X} \setminus \pi^{-1}(Z)).$$
But by [**18**, Corollary 1.4.5], we have that
$$X \setminus Z \cong \tilde{X} \setminus \pi^{-1}(Z). \tag{1}$$
Since $X$ is reduced, we must have that $\mathcal{O}_{\tilde{X}, \tilde{x}}$ is also reduced, for any point $\tilde{x} \notin \pi^{-1}(Z)$. In other words, that $\tilde{\mathcal{G}}(\tilde{X} \setminus \pi^{-1}(Z)) = 0$, as claimed.

To prove surjectivity, note that $V = \text{Im}(\pi)$ is a closed analytic subset of $X$, since $\pi$ is proper. On the other hand, by (1), it also contains the dense subset $X \setminus Z$, implying that $X = V$, since the former is reduced. ∎

**4.7. Definition.** Let $f \colon Y \to X$ be a map of affinoid varieties. We call $f$ *generically dense*, if the Zariski-closure $\text{cl}(f(\Sigma))$ of the image $f(\Sigma)$ of any irreducible component $\Sigma$ of $Y$, is an irreducible component of $X$. In terms of the corresponding morphism of affinoid algebras $A \to B$, this translates into any minimal prime of $B$ contracting to a minimal prime of $A$. For a map $f \colon Y \to X$ of quasi-compact rigid analytic varieties, we say that $f$ is *generically dense*, if there exist finite admissible affinoid coverings $\{Y_i\}_{i<n}$ of $Y$ and $\{X_i\}_{i<n}$ of $X$, such that $f(Y_i) \subset X_i$ and each of the maps $f_i = f|_{Y_i} \colon Y_i \to X_i$ is generically dense. For a point $x \in X$, we will say that $f$ is *locally generically dense* at $x$, if there exists an admissible affinoid $U$ of $X$ containing $x$, such that the restriction map $f^{-1}(U) \to U$ is generically dense.

**4.8. Lemma.** *Let $A \hookrightarrow B$ be an inclusion of affinoid algebras. Then $\dim A \leq \dim B$.*

*Proof.* Since $A$ is finite over a free Tate ring (necessarily of the same dimension), by Noether Normalization [**3**, 6.1.2. Corollary 2], we may already assume that $A = K\langle S \rangle$, where $S = (S_1, \ldots, S_d)$ are variables. Let $T = (T_1, \ldots, T_n)$ be more variables, such that $B = K\langle S, T \rangle / \mathfrak{n}$. We will prove by induction on $n$, that $\dim B \geq d$. There is nothing to prove if $\mathfrak{n} = 0$, so that in particular the case $n = 0$ is trivial. Hence let $f \in \mathfrak{n}$ be a non-zero element. By the Weierstrass Preparation Theorem, we can find a finite and injective map $K\langle S, T' \rangle \hookrightarrow K\langle S, T \rangle / (f)$, where $T' = (T_1, \ldots, T_{n-1})$. Let $\mathfrak{n}' = \mathfrak{n} \cap K\langle S, T' \rangle$ and let $B' = K\langle S, T' \rangle / \mathfrak{n}'$. We therefore have injective maps $K\langle S \rangle \hookrightarrow B' \hookrightarrow B$, where the latter map is finite, so that $\dim B' = \dim B$. We are now done by our induction hypothesis on $n$. ∎

**4.9. Lemma.** *Let $A \to C$ be a map of affinoid varieties turning $\operatorname{Sp} C$ into an affinoid subdomain of $\operatorname{Sp} A$. If $\mathfrak{P} \in \operatorname{Spec}(C)$ such that its contraction $\mathfrak{p} = \mathfrak{P} \cap A$ is a minimal prime of $A$, then also $\mathfrak{P}$ is a minimal prime of $C$.*

*Proof.* It is an easy exercise to see that it is enough to prove that $\mathfrak{P}$ is a minimal prime of $C/\mathfrak{p}C$. So we may already assume that $A$ is a domain and that $A \hookrightarrow C/\mathfrak{P}$. From **(4.8)** it follows that $\dim A \leq \dim C/\mathfrak{P}$. If $\mathfrak{P}$ were not minimal, then $\dim C/\mathfrak{P} < \dim C$. But since $\operatorname{Sp} C$ is an affinoid subdomain in $\operatorname{Sp} A$, we have that $\dim A = \dim C$, which would contradict the first inequality. Hence $\mathfrak{P}$ must be minimal, as we wanted to show. ■

**4.10. Lemma.** *Let $f: Y \to X$ be a map of quasi-compact rigid analytic varieties and let $U$ be an admissible open of $X$. If $f$ is generically dense, then also the restriction $f^{-1}(U) \to U$.*

*Proof.* Since the question is local, we may assume that $X = \operatorname{Sp} A$, $Y = \operatorname{Sp} B$ and $U = \operatorname{Sp} C$ are affinoid. Hence the affinoid algebra of $f^{-1}(U)$ is $C \widehat{\otimes}_A B$. Consider the following commutative diagram

$$\begin{array}{ccc} A & \xrightarrow{\varphi} & B \\ \delta \downarrow & & \downarrow \varepsilon \\ C & \xrightarrow{\varphi_0} & C \widehat{\otimes}_A B. \end{array}$$

Let $\mathfrak{Q}$ be a minimal prime of $C \widehat{\otimes}_A B$. Let $\mathfrak{P}$ (respectively, $\mathfrak{q}$) denote its contraction in $C$ (respectively, in $B$). Hence, since $\varepsilon$ is flat by [**3**, 7.3.2. Corollary 6], we have by [**10**, Theorem 15.1] that $\mathfrak{q}$ is a minimal prime of $B$. Since, by assumption $f$ is generically dense, we have that $\mathfrak{p} = \mathfrak{q} \cap A$ is then also minimal. But $\mathfrak{P} \cap A$ also equals $\mathfrak{p}$, so that $\mathfrak{P}$ must be minimal as well, by **(4.9)**. ■

**4.11. Proposition.** *Let $f: Y \to X$ be a map of quasi-compact rigid analytic varieties. Then $f$ is generically dense, if and only if, the inverse image of any nowhere dense locally closed analytic subvariety of $X$ is nowhere dense in $Y$. In particular, Zariski-open maps (i.e., maps sending Zariski-opens to Zariski-opens) are generically dense.*

*Proof.* $\Longrightarrow$. We claim that it is enough to prove this implication under the additional assumption that $X$ and $Y$ are affinoid. Indeed, for the general case, let $\{Y_i\}_i$ and $\{X_i\}_i$ be finite admissible affinoid coverings of $Y$ and $X$ respectively, such that each $f(Y_i) \subset X_i$ and $f_i = f|_{Y_i}$ is generically dense. Let $V$ be a nowhere dense locally closed analytic subvariety of $X$. Hence $V \cap X_i$ is nowhere dense in $X_i$ and hence by the affinoid case, $f_i^{-1}(V \cap X_i) = f^{-1}(V) \cap Y_i$ is nowhere dense in $Y_i$. Since the $Y_i$ form an admissible covering, it follows that then also $f^{-1}(V)$ is nowhere dense.

Therefore, we may assume that $X = \operatorname{Sp} A$ and $Y = \operatorname{Sp} B$ are affinoid. Let $Z$ be a nowhere dense locally closed analytic subvariety of $X$ defined over an affinoid subdomain $U = \operatorname{Sp} C$ of $X$ by an ideal $\mathfrak{a}$ of $C$. By **(4.10)**, $f^{-1}(U) \to U$ is also generically dense, so that without loss of generality, we may assume that $Z$ is even a closed analytic subvariety of $X$, i.e., $\mathfrak{a}$ is an ideal of $A$.

According to **(4.5)**, we have that $H_A(\mathfrak{a})$ is nilpotent and we need to show that the same holds for $H_B(\mathfrak{a}B)$. Suppose the contrary, so that there exists a minimal prime $\mathfrak{P}$ of $B$ not containing $H_B(\mathfrak{a}B)$. Let $n$ be such that $H_B(\mathfrak{a}B) = \operatorname{Ann}_B(\mathfrak{a}^n B)$. Hence we can find some $s \notin \mathfrak{P}$ such that $s\mathfrak{a}^n B = 0$. This implies that $\mathfrak{a}B \subset \mathfrak{P}$. Let $\mathfrak{p} = \mathfrak{P} \cap A$, so that then $\mathfrak{a} \subset \mathfrak{p}$. In particular, we have that

$$\operatorname{Ann}_A(\mathfrak{p}) \subset \operatorname{Ann}_A(\mathfrak{a}) \subset H_A(\mathfrak{a}) \subset \operatorname{nil}(A) \subset \mathfrak{p}$$

implying that $\mathfrak{p}^2 = 0$. But then also $\mathfrak{a}^2 = 0$, implying that $H_A(\mathfrak{a}) = A$, contradiction.

$\Longleftarrow$. We can always find finite admissible affinoid coverings $\{Y_i\}_i$ and $\{X_i\}_i$ of $Y$ and $X$ respectively, such that each $f(Y_i) \subset X_i$. We seek to show that the maps $f_i = f|_{Y_i}: Y_i \to X_i$ are generically dense. Therefore, let $\Sigma$ be an irreducible component of $Y_i$ and let $V = \operatorname{cl}_{X_i}(f(\Sigma))$, which necessarily has to be irreducible in $X_i$. Suppose that $V$ is not an irreducible component of $X_i$, so that by **(4.5)** $V$ must be nowhere dense in $X_i$. Since $V$ is then a nowhere dense locally closed analytic subset of $X$, it then follows, by our assumption, that $f^{-1}(V)$ is nowhere dense and whence also $f_i^{-1}(V)$. However, the latter set contains $\Sigma$ and hence is not nowhere dense by **(4.5)**, contradiction.

To prove the last statement, we only need to show that the inverse image $f^{-1}(Z)$ of a nowhere dense closed subset $Z$ of $X$ is nowhere dense in $Y$, if $f$ is Zariski-open. Suppose the contrary, so that

RIGID ANALYTIC FLATIFICATORS 17there exists a non-empty Zariski-open subset $U$ of $Y$ contained in $f^{-1}(Z)$. Hence $f(U) \subset Z$ and, by assumption, $f(U)$ is a non-empty Zariski-open subset of $X$, contradicting the nowhere density of $Z$.
∎

*Remark.* From the proof of $\Longleftarrow$, it follows that a map $f\colon Y \to X$ of quasi-compact rigid analytic varieties is nowhere dense, if and only if, for every pair of admissible affinoid coverings $\{Y_i\}_i$ and $\{X_i\}_i$ of $Y$ and $X$ respectively such that $f(Y_i) \subset X_i$, for each $i$, we have that each $f|_{Y_i}$ is generically dense. In other words, in the definition of generically density, we may replace the phrase *there exist finite admissible affinoid coverings such that ...* by the phrase *for every pair of admissible coverings such that ...* .

## 5. Fibre Lemma

**5.1. Lemma.** *Let $S$ be a Noetherian ring, $\mathfrak{a}$ a proper ideal in $S$ for which the $\mathfrak{a}$-adic topology on $S$ is separated and let $S \to A$ be a morphism of rings. Suppose that $A$ has a founded and flat $\mathfrak{a}$-representation*
$$0 \to H \longrightarrow F \xrightarrow{\pi} A \to 0$$
*By foundedness, there exists a smallest ideal $\mathfrak{f}$ for which $H \subset \mathfrak{f}F$. Suppose that $\mathfrak{f}$ is invertible, i.e., $\mathfrak{f} = (x)$ with $x \in S$ a non zero-divisor. Then there exists an $h \in A \setminus \mathfrak{a}A$, such that $xh = 0$ in $A$.*

*Proof.* We must have that
$$H \not\subset (x)\mathfrak{a}F,$$
since otherwise $\mathfrak{f} = (x) \subset (x)\mathfrak{a}$, which would imply that $(x) \subset \mathfrak{a}^n$, for every $n$ and hence that $x = 0$ by the separatedness condition. So, we can find an element $k \in H$, with $k = xf$ and $f \notin \mathfrak{a}F$. Since $\pi$ is an $\mathfrak{a}$-representation, the latter implies that $\pi(f) \notin \mathfrak{a}A$. Hence if we set $h = \pi(f)$, we get the desired element, since
$$0 = \pi(k) = xh.$$
∎

**5.2. Theorem (Fibre lemma).** *Let $f\colon Y \to X$ be a map of affinoid varieties and let $x$ be a point of $X$. Suppose that $f$ has a globally defined universal rigid analytic flatificator $Z$ at $x$ (i.e., $Z$ is a closed analytic subvariety of $X$). Let $\pi\colon \tilde{X} \to X$ denote the blowing-up map with centre $Z$ and let*

$$\begin{array}{ccc} \tilde{Y} & \xrightarrow{\theta} & Y \\ \tilde{f}\downarrow & & \downarrow f \\ \tilde{X} & \xrightarrow{\pi} & X \end{array}$$

*be the commutative diagram of the strict transform.*

*Then we have, for each point $\tilde{x} \in \pi^{-1}(x)$, a strict closed immersion (i.e., a closed immersion which is not an isomorphism)*

(1) $$\tilde{f}^{-1}(\tilde{x}) \subsetneq f^{-1}(x)$$

*of closed analytic subvarieties.*

*Remark.* Note that there is always a closed immersion
$$\tilde{f}^{-1}(\tilde{x}) \subset f^{-1}(x).$$

Indeed, $\tilde{Y}$, by definition of strict transform, is the blowing up of $Y$ with centre $f^{-1}(Z)$. By [**18**, Proposition 3.1.2], we have a closed immersion

(2) $$i\colon \tilde{Y} \hookrightarrow Y \times_X \tilde{X},$$

which induces a closed immersion of fibres
$$\tilde{f}^{-1}(\tilde{x}) \hookrightarrow \gamma^{-1}(\tilde{x}) \cong f^{-1}(x),$$

where $\gamma\colon Y\times_X \tilde X \to \tilde X$ denotes the canonical projection. See also below in the proof, for some more explicit calculations. We emphasize that we will prove (1), for the fibres considered as subvarieties rather than as subsets (i.e., it may be that the underlying closed analytic sets are equal but they differ in their analytic structure).

*Proof.* Let $A$ (respectively, $B$) be the affinoid algebras of $X$ (respectively, $Y$). Let $\mathfrak{m}$ be the maximal ideal of $A$ corresponding to the point $x$. Then the fibre of $x$ over $f$ is given as $\operatorname{Sp}(B/\mathfrak{m}B)$. Let $\mathfrak{f}$ be the ideal defining $Z$. Let $\tilde U = \operatorname{Sp} \tilde C$ be an affinoid admissible open of $\tilde X$ containing $\tilde x$. From (2), we then obtain that there is a closed immersion

$$\tilde f^{-1}(\tilde U) \hookrightarrow Y\times_X \tilde U. \tag{3}$$

The latter space is affinoid, with affinoid algebra $D = B\widehat{\otimes}_A \tilde C$. Let $I$ be the ideal in $D$ giving this closed immersion (3), so that $\tilde f^{-1}(\tilde U)$ is affinoid with affinoid algebra $\tilde D = D/I$. Let us calculate the various fibres. Let $\tilde{\mathfrak{m}}$ be the maximal ideal in $\tilde C$ defining $\tilde x$. Hence, via the map $\tilde f^{-1}(\tilde U) = \operatorname{Sp}\tilde D \to \tilde U = \operatorname{Sp}\tilde C$, we see that

$$\tilde f^{-1}(\tilde x) = \operatorname{Sp}(\tilde D/\tilde{\mathfrak{m}}\tilde D).$$

On the other hand, we have that

$$\gamma^{-1}(\tilde x) = \operatorname{Sp}(D/\tilde{\mathfrak{m}}D).$$

But

$$D/\tilde{\mathfrak{m}}D \cong (B/\mathfrak{m}B)\widehat{\otimes}_{A/\mathfrak{m}}(\tilde C/\tilde{\mathfrak{m}}) \cong B/\mathfrak{m}B$$

where the latter algebra is the affinoid algebra of $f^{-1}(x)$. This proves that $\gamma^{-1}(\tilde x) \cong f^{-1}(x)$ and that $\tilde f^{-1}(\tilde x)$ is a closed subvariety of $f^{-1}(x)$ given by the (image of the) ideal $I$ in $D/\tilde{\mathfrak{m}}D$. So we proved (1), if we can show that $I \neq 0$ in $D/\tilde{\mathfrak{m}}D$.

By definition of blowing up, we have that the ideal $\mathfrak{f}\mathcal{O}_{\tilde X}$ has become invertible. By shrinking $\tilde U$ if necessary, we may already assume that $\mathfrak{f}\tilde C$ is generated by a non zero-divisor of $\tilde C$. Moreover, since $\mathfrak{f}$ is universal, we have that $\mathfrak{f}\tilde C$ is a flatificator of the base change morphism $\tilde C \to D = B\widehat{\otimes}_A\tilde C$. Hence from **(5.1)** it follows that there exists an $h \in D \setminus \tilde{\mathfrak{m}}D$, for which $h\mathfrak{f} = 0$ in $D$, since $\mathfrak{f} \subset \tilde{\mathfrak{m}}$ and since we always may assume that $\tilde C \to D$ has a founded and flat $\tilde{\mathfrak{m}}$-representation by **(3.3)**, after perhaps some further shrinking of $\tilde U$. By the definition of strict transform, also $\mathfrak{f}\tilde D$ is invertible. Since $h\mathfrak{f} = 0$ in $D$ whence in $\tilde D$, we conclude that $h = 0$ in $\tilde D$ and therefore $h \in I$. But $h \neq 0$ in $D/\tilde{\mathfrak{m}}D$, proving that $I \neq 0$ in $D/\tilde{\mathfrak{m}}D$. ∎

**5.3. Theorem.** *Let $f\colon Y \to X$ be a map of quasi-compact rigid analytic varieties and let $x \in X$ be an arbitrary point. Suppose that $X$ is reduced. If $f$ is not flat at $x$, then there exists a local blowing up $\pi\colon \tilde X \to X$ with nowhere dense centre $F$ containing $x$, such that $f^{-1}(F) \to F$ is flat and such that, if*

$$\begin{array}{ccc} \tilde Y & \xrightarrow{\theta} & Y \\ \tilde f \downarrow & & \downarrow f \\ \tilde X & \xrightarrow{\pi} & X \end{array}$$

*denotes the diagram of the strict transform of $f$ under $\pi$, we have, for every $\tilde x \in \pi^{-1}(x)$, a non-trivial closed immersion of fibres*

$$\tilde f^{-1}(\tilde x) \subsetneq f^{-1}(x). \tag{1}$$

*Proof.* Since the question is local, we may assume that $X$ and $Y$ are affinoid. According to **(3.7)**, we can find an affinoid subdomain $U$ of $X$ containing $x$, such that the restriction map $f^{-1}(U) \to U$ has a universal rigid analytic flatificator $Z$ at $x$, which we can assume to be globally defined. Apply **(5.2)** to this situation to obtain the required strict inclusion of fibres. Note also that by definition of a flatificator the condition that $f^{-1}(Z) \to Z$ be flat, is automatically met (and hence also for any closed analytic subvariety of $Z$).

Therefore, the only thing that remains to be shown, is that we can replace $Z$ by a smaller centre $F$ which is nowhere dense without violating (1). Without loss of generality, we can always, after shrinking $U$ if necessary, assume that $x$ lies on each irreducible component of $U$. Suppose that $V$ is an irreducible component of $U$ contained in $|Z|$. Let $X_0$ denote the analytic complement $\complement_U(V)$ and $Z_0 = Z \cap X_0$. Note that since $|Z| \subsetneq U$, where the latter is reduced, and since $x$ lies on all irreducible components of $U$, we must have that $x \in X_0$. Let $\pi_0 \colon \tilde{X}_0 \to U$ be the blowing up of $U$ with centre $Z_0$. By **(4.4)**, we have an isomorphism

$$\pi^{-1}(X_0) \cong \pi_0^{-1}(X_0).$$

In particular, since $x \in X_0$, we see that the fibre of $x$ under $\pi$ is isomorphic with the fibre of $x$ under $\pi_0$.

Let us now look at the strict transform of $Y$ under $\pi_0$, i.e., the blowing up of $f^{-1}(U)$ with centre $f^{-1}(Z_0)$. Put $Y_0 = \complement_{f^{-1}(U)}(f^{-1}(V))$. Then we have closed immersions[2]

$$f^{-1}(Z) \cap Y_0 \hookrightarrow f^{-1}(Z_0) \hookrightarrow f^{-1}(Z).$$

Again, by applying **(4.4)** and the same reasoning as above, we see that the fibre above any $\tilde{x} \in \pi^{-1}(x)$ has not changed. In other words, the closed immersion (1) of the fibres remains the same when we replace the local blowing up $\pi$ by the local blowing up $\pi_0$. Proceeding in this way, we can eliminate all irreducible components of $U$ that are lying in the centre of blowing up, to obtain a centre which is nowhere dense by **(4.5)**. ∎

*Remark.* We can even be more precise in the above statement. Suppose $y \in Y$ is a point lying above $x$ in which $f$ is not flat. Then there exists an admissible affinoid $U \subset Y$ containing $y$, such that

$$\tilde{f}^{-1}(\tilde{x}) \cap \theta^{-1}(U) \subsetneq f^{-1}(x) \cap U.$$

Indeed, just observe that the flatificator of the restriction $f|_U \colon U \to X$ at $x$ is contained in the flatificator of $f$ at $x$ and apply now the theorem to this new map $f|_U$.

**5.4. Theorem.** *Let $f \colon Y \to X$ be a map of rigid analytic varieties with $X$ reduced and let $x \in X$. If the following two conditions hold*

(1) *$f$ is locally generically dense at $x$,*
(2) *The fibre $f^{-1}(x)$ is a single point $\{y\}$ with its reduced closed subspace structure,*

*then $f$ is flat at $x$ (or, equivalently, flat in $y$).*

*Proof.* Suppose that $f$ is not flat at $x$. Using **(5.3)**, we can find an admissible affinoid $U$ around $x$ and a nowhere dense closed analytic subvariety $Z$ of $U$ such that the following holds. Let $\pi \colon \tilde{X} \to U$ denote the blowing up of $U$ with centre $Z$. Let

(1)
$$\begin{array}{ccc} \tilde{Y} & \xrightarrow{\theta} & f^{-1}(U) \\ \tilde{f} \downarrow & & \downarrow f \\ \tilde{X} & \xrightarrow{\pi} & U \end{array}$$

be the diagram of the strict transform. Then, for each point $\tilde{x} \in \pi^{-1}(x)$, the fibre $\tilde{f}^{-1}(\tilde{x})$ is strictly contained in $f^{-1}(x) = \{y\}$, where the latter has the reduced subspace structure. In other words, we must have that each of these fibres $\tilde{f}^{-1}(\tilde{x})$ is empty. Therefore, we found that

(2) $$\pi^{-1}(x) \cap \tilde{f}(\tilde{Y}) = \emptyset.$$

By **(4.10)** and our first assumption, we can even assume that the map $f^{-1}(U) \to U$ is generically dense. Using **(4.11)** we obtain that then $f^{-1}(Z)$ is nowhere dense in $f^{-1}(U)$. By **(4.6)**, $\theta$, being the blowing up of $f^{-1}(U)$ with centre $f^{-1}(Z)$, is surjective, so that there exists a $\tilde{y} \in \tilde{Y}$, such that $\theta(\tilde{y}) = y$. This implies, by the commutativity of (1), that $\tilde{f}(\tilde{y})$ lies in the fibre $\pi^{-1}(x)$. But this clearly contradicts (2), so that $f$ must be flat at $x$. ∎

---

[2]We have used here the following preservation property of analytic complements under morphisms: let $f \colon Y \to X$ be a map of rigid analytic varieties and let $V$ be a closed analytic subvariety of $X$, then

$$f(\complement_Y(f^{-1}(V))) \subset \complement_X(V).$$

**5.5. Corollary.** *Let $\varphi\colon A \hookrightarrow B$ be an injective morphism of affinoid domains. Let $\mathfrak{m}$ be a maximal ideal of $A$. If $\mathfrak{m}B$ is a maximal ideal in $B$, then $B$ is flat over $A$ at $\mathfrak{m}$ (meaning that $A_{\mathfrak{m}} \to B_{\mathfrak{m}B}$ is flat).*

*Proof.* We only have to verify the conditions of **(5.4)**. Since both algebras are domains and since $\varphi$ is injective, $\varphi$ is generically dense at any point, by **(4.10)**. The fibre at $\mathfrak{m}$ equals $\mathrm{Sp}(B/\mathfrak{m}B) = \mathrm{Sp}\,K$, so that (2) of **(5.4)** also holds. ∎

## Appendix: Flatificators in Berkovich Spaces

**A.1. Definition.** We will not give a complete account on Berkovich spaces, but just give a brief outline, so that the reader will convince himself that a theory of flatificators, analogous to the above, can be developed in this new setup. For a more detailed description, we have to refer to Gardener's D.Phil. thesis [**5**]. For an introduction to the general theory of Berkovich spaces, see also Berkovich's own work [**1**] and [**2**] or the more accessible renderings by Schneider [**14**] and Schneider and van der Put [**15**].

Let $X = \mathrm{Sp}\,A$ be an affinoid variety. An *analytic point* $x$ of $X$ is a continuous $K$-algebra morphism $x\colon A \to F$ to a complete extension field $F$ of $K$ (i.e., such that $F$ comes with a complete norm, extending the one on $K$). For convenience sake, we will require also that $x(A)$ generates a dense subfield of $F$.[3] An ordinary point of $X$ corresponds to a maximal ideal of $A$ and hence induces a morphism $A \to K$, thus becoming an analytic point. These will be called *geometric points*. An affinoid subdomain $U = \mathrm{Sp}\,C$ is called an *affinoid neighbourhood* of an analytic point $x$, if the morphism $x\colon A \to F$ factors through a morphism $C \to F$, still denoted $x$ (which is then necessarily unique). Two analytic points are said to be congruent if they have the same system of affinoid neighbourhoods. The *K-affinoid Berkovich space* associated to $X$ (or to $A$) is defined to be the set $\mathbb{M}(X)$ of all congruence classes of analytic points of $X$. The latter set is topologized by putting the weakest topology on it rendering each map

$$\tilde{f}\colon \mathbb{M}(X) \to \mathbb{R}\colon x \mapsto |x(f)|$$

continuous, where $f \in A$ and $x\colon A \to F$ (the norm is the one in $F$). This turns $\mathbb{M}(X)$ into a compact space and $X$ embeds in $\mathbb{M}(X)$ as an everywhere dense subset. Finally, to define the $K$-analytic structure, we make $\mathbb{M}(X)$ into a locally ringed $K$-space by defining a structure sheaf $\mathcal{O}_{\mathbb{M}(X)}$.

If $U \subset X$ is an affinoid neighbourhood of an analytic point $x$ of $X$, then $\mathbb{M}(U)$ is the closure of $U$ in $\mathbb{M}(X)$. It is called a *wide* affinoid neighbourhood of $x$, if $\mathbb{M}(U)$ contains an open neighbourhood of $x$ in $\mathbb{M}(X)$. If $x$ is geometric, any affinoid subdomain is wide. The compactness of the Berkovich space translates into the following property of the associated affinoid variety (see [**2**, Lemma 1.6.2]): if $\{U_i\}_i$ is a collection of affinoid subdomains of $X$, such that, for each analytic point $x$ of $X$, one of the $U_i$ is a wide affinoid neighbourhood of $x$, then the covering is admissible (in the Grothendieck topology on $X$) and whence already finitely many $U_i$ cover $X$. Beware that the above condition is stronger than just requiring $\{U_i\}_i$ to be a covering, since the latter only means that the geometric points are 'covered'.

There exist global versions of the Berkovich space defined for general rigid analytic varieties. To our purposes, we will only be concerned with separated quasi-compact rigid analytic varieties. To any quasi-compact rigid analytic variety $X$, one can associate a *K-analytic Berkovich space* $\mathbb{M}(X)$ (the latter space is then locally compact, paracompact and Hausdorff). In this paper a Berkovich space will always be of the latter form, contrary to the more general usage of the terminology in the aforesaid sources, where the interested reader can find the details. Let us just mention the following facts: any map $f\colon Y \to X$ between quasi-compact rigid analytic varieties uniquely determines a map $\mathbb{M}(f)\colon \mathbb{M}(Y) \to \mathbb{M}(X)$ between the corresponding Berkovich spaces. The former is flat if the latter is, where $\mathbb{M}(f)$ is called *flat* at a point $y \in \mathbb{M}(Y)$ if the map of local rings $\mathcal{O}_{\mathbb{M}(X),\mathbb{M}(f)(y)} \to \mathcal{O}_{\mathbb{M}(Y),y}$ is flat.

The definition of flatificator and universal flatificator are now similar. Let $f\colon \mathbb{Y} \to \mathbb{X}$ be a map of Berkovich spaces, let $x$ be a point of $\mathbb{X}$ and let $\mathbb{Z}$ be a locally closed analytic subspace of $\mathbb{X}$ containing $x$. Then we say that $\mathbb{Z}$ is the *Berkovich flatificator* of $f$ at $x$, if $f^{-1}(\mathbb{Z}) \to \mathbb{Z}$ is flat, and for every locally closed analytic subspace $\mathbb{V}$ of $\mathbb{X}$ containing $x$ such that $f^{-1}(\mathbb{V}) \to \mathbb{V}$ is flat at $x$, we have that $\mathbb{V} \subset \mathbb{Z}$ at $x$. It is called a *universal Berkovich flatificator* of $f$ at $x$, if it remains the flatificator at all points of the fibre of $x$ over any base change. As in the rigid analytic case, for a map to have a Berkovich flatificator is equivalent with the corresponding morphism $\mathcal{O}_{\mathbb{X},x} \to (f_*\mathcal{O}_{\mathbb{Y}})_x$ having a flatificator $\mathfrak{A}$ such that $\mathcal{O}_\mathbb{X}(\mathbb{U})/\mathfrak{A} \to \mathcal{O}_\mathbb{Y}(f^{-1}(\mathbb{U}))/\mathfrak{A}\mathcal{O}_\mathbb{Y}(f^{-1}(\mathbb{U}))$ is flat, where $\mathbb{U}$ is an open neighbourhood of $x$ on which $\mathfrak{A}$ is defined. We can now formulate the Berkovich analogue of our Theorem **(3.7)**.

---

[3]The reader should convince himself that this is a harmless extra condition, because we will only be interested in congruence classes of analytic points, as to be defined below.

**A.2. Theorem.** *Let $f \colon \mathbb{X} \to \mathbb{Y}$ be a map of Berkovich spaces and let $x$ be a point of $\mathbb{X}$. Then there exists a wide affinoid neighbourhood $\mathbb{U}$ of $x$, such that the restriction $f^{-1}(\mathbb{U}) \to \mathbb{U}$ has a universal Berkovich flatificator.*

*Proof ( [5, Proposition 5.8] ).* First we need to prove the Berkovich-analogue of **(3.3)**, i.e., replacing, in the statement of , the geometric point $x$ by an arbitrary analytic point and the affinoid subdomain $U$ containing $x$ by a wide affinoid neighbourhood of $x$. Note that an analytic point $x$ is a morphism $A \to F$, where $F$ is a complete extension field of $K$ and hence its kernel is a prime ideal $\mathfrak{p}$ of $A$. Therefore, replace in the proof of , the maximal ideal $\mathfrak{m}$ by the prime ideal $\mathfrak{p}$. Similarly, the affinoid subdomain $U$ defined by (2) has to be modified by replacing the first condition in (2) by $|a_\nu(z)| \geq \alpha - \varepsilon/2$ (and a similar adaptation for the case $n = 0$), thus ensuring that it becomes a wide affinoid neighbourhood of $x$. One should also extend the definition of a strongly $\mathcal{O}_X$-module $\mathcal{F}$ to all (geometric or analytic) fibres being founded and observe that **(1.2)** remains valid. The rest of the proof can now be copied verbatim to entail the Berkovich-analogue of **(3.3)**. To complete the proof of the theorem, we can now use the same proof as for **(3.7)**, since nowhere we have exploited that the ideal corresponding to a point is maximal, but only that it is prime. ∎

**A.3. Corollary.** *Flatness for Berkovich maps is an open condition in the source.*

*Proof ( [5, Proposition 5.7] ).* Same proof as for **(3.8)**. ∎

In order to formulate the Berkovich variant of **(5.3)**, we need to say a bit more about fibres.

**A.4. Definition.** Let $A \to B$ be a map of $K$-affinoid algebras, let $\mathbb{X} = \mathbb{M}(\operatorname{Sp} A)$ and $\mathbb{Y} = \mathbb{M}(\operatorname{Sp} B)$ be the corresponding $K$-affinoid Berkovich spaces and let $f \colon \mathbb{Y} \to \mathbb{X}$ be the corresponding map. Let $x \colon A \to F$ be a point of $\mathbb{X}$, i.e., an analytic point of $X = \operatorname{Sp} A$. Let $y \in f^{-1}(x)$ with $y \colon B \to G$ a point of $\mathbb{Y}$, then the additional requirement that $x(A)$ generates a dense subfield of $F$, ensures that there is a continuous embedding $F \hookrightarrow G$. Therefore, $y$ can be viewed as an analytic point of $Y_x = \operatorname{Sp}(B \widehat{\otimes}_A F)$ and, conversely, any analytic point of $Y_x$ induces an analytic point of $Y = \operatorname{Sp} B$. We will call the $F$-affinoid Berkovich space $\mathbb{M}(Y_x)$ the *fibre* over $x$ and we denote it as $\mathbb{Y}_x$. Note that $\mathbb{Y}_x$ is naturally homeomorphic with $f^{-1}(x)$. However, in general, it will no longer be $K$-analytic anymore, but $F$-analytic, for some complete extension field $F$ of $K$ (although for geometric points $K = F$, so that in that case the above definition agrees with the one used in Section 5). It is this fact that has to be taken into account when considering the Fibre Lemma. We have the following easy to prove result.

**A.5. Fact.** *Let $f \colon \mathbb{Y} \to \mathbb{X}$ and $g \colon \mathbb{Z} \to \mathbb{X}$ be maps of $K$-affinoid Berkovich spaces. Let $x \colon A \to F$ (respectively, $z \colon C \to G$) be a point of $\mathbb{X}$ (respectively, of $\mathbb{Z}$), where $A$ (respectively, $C$) is the affinoid algebra associated to $\mathbb{X}$ (respectively, to $\mathbb{Z}$). If $g(z) = x$ and if $f'$ denotes the base change $\mathbb{W} = \mathbb{Z} \times_\mathbb{X} \mathbb{Y} \to \mathbb{Z}$ of $f$, then we have an isomorphism*

$$\mathbb{W}_z \cong \mathbb{Y}_x \times_F G$$

*of $G$-analytic Berkovich spaces (where the above fibre product just means extension of scalars).*

It will come as no surprise that one can define blowing up maps and local blowing up maps in the category of Berkovich spaces as well and that the functor $\mathbb{M}(\cdot)$ maps rigid analytic blowing up maps to Berkovich blowing up maps, whenever this makes sense. We finally come to the last Berkovich variant, where **(5.3)** gets translated into the following theorem.

**A.6. Theorem.** *Let $f\colon \mathbb{Y} \to \mathbb{X}$ be a map of $K$-analytic Berkovich spaces and let $x\colon A \to F$ be an arbitrary point of $\mathbb{X}$. Suppose that $\mathbb{X}$ is reduced and take some compact subset $\mathbb{L}$ of the fibre $f^{-1}(x)$. If $f$ is not flat in each point of $\mathbb{L}$, then there exists a local blowing up $\pi\colon \tilde{\mathbb{X}} \to \mathbb{X}$ with nowhere dense centre $\mathbb{Z}$ containing $x$, such that $f^{-1}(\mathbb{Z}) \to \mathbb{Z}$ is flat and such that, if*

$$\begin{array}{ccc} \tilde{\mathbb{Y}} & \xrightarrow{\theta} & \mathbb{Y} \\ \tilde{f}\downarrow & & \downarrow f \\ \tilde{\mathbb{X}} & \xrightarrow{\pi} & \mathbb{X} \end{array}$$

*denotes the diagram of the strict transform of $f$ under $\pi$, we have, for every $\tilde{x}\colon \tilde{A} \to \tilde{F}$ in $\tilde{\mathbb{X}}_x$ (where $\mathrm{Sp}(\tilde{A})$ is some wide affinoid neighbourhood of $\tilde{x}$), a closed immersion of fibres*

$$\tilde{\mathbb{Y}}_{\tilde{x}} \hookrightarrow \mathbb{Y}_x \times_F \tilde{F} = (\tilde{\mathbb{X}} \times_{\mathbb{X}} \mathbb{Y})_{\tilde{x}}$$

*which is strict at some point in $\{\tilde{x}\} \times_{\mathbb{X}} \mathbb{L}$. In fact, this inclusion remains strict after extension of scalars to an arbitrary complete extension field $G$ of $\tilde{F}$.*

*Proof.* Apart from extending the scalars, we can copy the rigid analytic Fibre Lemma **(5.2)** verbatim to the Berkovich setup (just replace the maximal ideals corresponding to geometric points by the prime ideals corresponding to analytic points). Likewise for **(5.3)**, since the whole machinery on analytic complements in Section 4 can be transformed to its Berkovich equivalent by applying the fully faithful functor $\mathbb{M}(\cdot)$. (For instance, the *analytic Berkovich complement* of a closed analytic subspace $\mathbb{F}$ of the $K$-affinoid Berkovich space $\mathbb{X}$ is defined as $\mathbb{M}(\complement_X(F))$, where $F$ and $X$ are the affinoid varieties associated to $\mathbb{F}$ and $\mathbb{X}$, etc.). Finally, use the remark after **(5.3)** to show that the strict inclusion (1) already holds in a point above $\mathbb{L}$.　∎

FIELDS INSTITUTE
222 COLLEGE STREET
TORONTO, ONTARIO, M5T 3J1 (CANADA)